%% file: MU.tex
\documentclass{amsart}

\include{KU-0header}

\begin{document}


\title[RBMs for Parabolic PDEs with Parameter Functions]%
    {A Reduced Basis Method for Parabolic Partial Differential Equations with Parameter Functions and Application to Option Pricing}


\author{Antonia Mayerhofer}

\address{
Antonia Mayerhofer\\
University of Ulm\\
Graduiertenkolleg 1100\\
Helm\-holtz\-strasse {22}\\
D-89069 Ulm, Germany}

\email{antonia.mayerhofer@uni-ulm.de}


\author{Karsten Urban}

\address{
Karsten Urban\\
University of Ulm\\
Institute for Numerical Mathematics\\
Helm\-holtz\-strasse {20}\\
D-89069 Ulm, Germany}

\email{karsten.urban@uni-ulm.de}

\begin{abstract}
 We consider the Heston model as an example of a parameterized parabolic partial differential equation. A space-time variational formulation is derived that allows for parameters in the coefficients (for calibration) as well as choosing the initial condition (for option pricing) as a parameter function. A corresponding discretization in space and time amd initial condition is introduced and shown to be stable. Finally, a Reduced Basis Method (RBM) is introduced that is able to use parameter functions also for the initial condition. Corresponding numerical results are shown.
\end{abstract}

\keywords{Option pricing, parabolic problems, reduced basis method, error estimates}

\subjclass[2010]{%
    35K85,
    49J40,
    65M15,
    91B25, 
    91G80
    }

\thanks{The authors have been supported by the Deutsche Forschungsgemeinschaft (DFG) under GrK 1100. We are grateful to Robert Stelzer for helpful discussions and comments.}

\date{\today}

\maketitle


\mysection{Introduction}
\label{Sec:0}

Calibration and pricing are two standard tasks in numerical finance. Using a PDE-model based upon the Feynman-Kac theorem (see also Proposition \ref{Prop:6.1} below), calibration amounts fitting unknown parameters in the partial differential equation (PDE) to historical market data in a least-squares sense e.g.\ by some optimization scheme. This usually requires solving the same PDE with many different values of the parameters (given by the iterative optimization method). Such a problem is also termed `multi-query'. 

We consider the Heston model, which allows for a non-constant volatility, \cite{Heston}. Having in mind that every model is wrong to a certain extend, constant volatilities are of course a severe limitation, which has been our motivation for considering the Heston model here. The Heston model leads to a system of parabolic diffusion-convection-reaction PDEs with non-constant coefficients, which makes the model also interesting from a numerical point of view. 

Moreover, as we shall see below, the Heston model already contains 4 parameters, i.e., calibration requires solving a 4-dimensional PDE-constrained optimization problem. Of course, by extending the model, the number of parameters can easily be enlarged leading to truly high-dimensional parameter spaces.

In the PDE-framework, pricing amounts solving the same kind of parabolic problem with different initial conditions determined by different payoff functions. The payoff is usually determined by the specific design of the financial product, so that a flexibility in the choice of the payoff function enables efficient product design.

Both problems, calibration and pricing need to be solved rapidly, often `online' e.g.\  in `real time' as for risk management.

Such kind of rapid online multi-query problems are in the scope of Reduced Basis Methods (RBM), that have been developed and widely investigated for parametric PDEs, in particular of elliptic and parabolic type, see, e.g.\ \cite{HaasdonkLuminy,RP} and references therein. Roughly speaking, the RBM involves an offline and an online computing phase. In the offline phase, a detailed, but costly numerical model (sometimes called `truth'), is used to construct a reduced space. The construction is based upon a reliable and effective error estimate which in addition can be evaluated efficiently. Moreover, this error estimate allows to compute an a posteriori bound, which makes RB-approximations certified.

Similar to \cite{CRAS,UP2}, we suggest a space-time variational formulation of the parabolic PDE-system. This leads to a Petrov-Galerkin method containing space and time. Obviously, the time is added as additional dimension, but it allows for sharp error bounds for a corresponding discretization in space and time. Moreover, similar to \cite{CRAS,UP2} we show that a specific choice of trial and test spaces in the Petrov-Galerkin method lead to a Crank-Nicholson scheme also for the case of varying initial conditions. This requires also a discretization of the initial condition so that we end up with a detailed discrete problem in space and time and intitial condition, which is used as `truth' for a RBM. It turns out that the space-time approach transfers an initial condition parameter function into a parameter within the PDE-coefficients.  

Even though RBM have already been used in numerical finance, see e.g.\ \cite{CLP,HSW,SaS}, and even though the current paper is clearly only a first step, we think to the best of our knowledge that we introduce some novel aspects: (a) We use the initial condition of an evolution problem as a parameter within a RBM, which has not been considered so far; (b) The initial condition is a parameter \emph{function}, which in particular means that the parameter space is infinite-dimensional. We solve this problem by expanding the initial condition in an appropriate basis, which is somehow similar to e.g.\ \cite{HUW,HS}. (c) A space-time variational formulation is used to treat a non-homogeneous initial condition by a two-step method, namely first approximating the intial condition and then using this as a parameter in the right-hand side for determining the evolution. We note, that this is \emph{not} restricted to linear problems only, see Remark \ref{Rem:nonlin} and \S \ref{Sec:nonlinear} below.

The remainder of this paper is organized as follows. In Section \ref{Sec:1}, we collect preliminaries on the Heston model and Reduced Basis Methods (RBM). Section \ref{Sec:3} is devoted to the review of the space-time variational formulation, in particular for parabolic parametric partial differential equations. This includes the error estimate from \cite{CRAS,UP2}, which is the basis for the RB error estimate. Next, we introduce a discretization. In space and time, we follow \cite{CRAS,UP2} and use finite elements. In addition, we introduce a discretization for the initial condition and show that we can decompose the numerical computations into (1) a Galerkin-type problem for the approximation of the initial condition and (2) a time-stepping scheme for the evolution.

The next step is the introduction of a Reduced Basis Method (RBM) in particular for parameter functions in Section \ref{Sec:4}. Again, we derive a two-step method, here for the construction of the reduced basis. We derive stable RB trial and test spaces for the reduced Petrov-Galerkin scheme and show an error estimate. We present numerical results in Section \ref{Sec:5} and close with a summary and outlook in Section \ref{Sec:6}.

\mysection{Preliminaries}
\label{Sec:1}

\subsection{Calibration and pricing within the Heston Model}
\label{Sec:1.1}
The Heston model is a well-known model for option pricing allowing also a non-constant volatility of the underlying. It was invented 1993 by Heston, \cite{Heston}. The  Stochastic Differential Equations (SDEs) for the asset price $S_{t}$ and the volatility $\nu_{t}$ are assumed to be
\begin{align}
    dS_{t} &= r  S_{t} dt +  \sqrt{\nu_{t}} S_{t} dz_1(t),
    & 
    d\nu_{t} &= \kappa [\theta - \nu_{t}] dt + \sigma \sqrt{\nu_{t}} d z_2(t),
    \label{Eq:SDEHeston}
\end{align}
where $z_{1}$, $z_{2}$ are Wiener processes with correlation $\rho$, $r$ is the return rate of the asset, $\kappa$ the mean reversion rate to the long term variance $\theta$ and $\sigma$ is the volatility of the volatility. In particular, the instantaneous variance $\nu_{t}$ is modeled as a CIR (Cox-Ingersoll-Ross) process, \cite{CIR}. Finally, the model implies that  $z_1 = \sqrt{1-\rho^2} dz_3 + \rho dz_2$ with  independent Brownian motions $z_2$, $z_3$.

The following partial differential equation (PDE) representation is well-known.

\begin{prop}[Feynman-Kac formula, \cite{Oe}]\label{Prop:6.1}
Let $X_t$ be an $n$-dimensional It\^{o} diffusion, i.e.,  $d X_t = b(X_t) dt + \sigma(X_t) dB_t$, $f \in C^{2}(\mathbb{R}^n)$ be a compactly supported payoff and $q \in C(\mathbb{R}^n)$ be lower bounded. For  $u(t,x) := \mathbb{E} [e^{-\int_0^t q(X_s) ds} f(X_t)| X_t=x]$ it holds for $x\in \mathbb{R}^n$ that $u(0,x)= f(x)$ and for all $t\in (0,T]$, $x\in \mathbb{R}^n$ that
\begin{align*}
    \frac{\partial u(t,x)}{\partial t} 
        &= \frac{1}{2} \sum_{i,j=1}^n (\sigma(x) \sigma^T (x))_{i,j} \frac{\partial^2 u(t,x)}{\partial x_i \partial x_j} + \sum_{i=1}^n b_i(x) \frac{\partial u(t,x)}{\partial x_i} - q(x)\, u(t,x).
        \hfill\qed
\end{align*}
\end{prop}

For applying the Feynman-Kac formula to the Heston model, we consider the process $X_t:= (y_t, \nu_t)^T = (\log(S_t), \nu_t)^T$, $n=2$. Then, we can rewrite the two SDEs in \eqref{Eq:SDEHeston} as a system
\begin{align*}
    dX_t = 
    d \begin{pmatrix} y_t \\ \nu_t \end{pmatrix} 
    &= \begin{pmatrix} r 
    - \frac{1}{2} \nu_{t} \\ \kappa [\theta - \nu_{t}] \end{pmatrix} dt 
    + \begin{pmatrix} \sqrt{\nu_{t}} \sqrt{1-\rho^2} & \sqrt{\nu_{t}} \rho \\ 
    0 & \sigma \sqrt{\nu_{t}} \end{pmatrix} \begin{pmatrix} dz_3 \\ dz_2 \end{pmatrix} \\
    &=: b(X_t)\, dt + \sigma(X_t)\, dB_t.
\end{align*}
Using Proposition \ref{Prop:6.1} then yields the PDE (in $x$) for the Heston model of the form
\begin{equation}
    \label{div}
    \frac{\partial u(t)}{\partial t} - \Div( \uuline\alpha(t)  \nabla u(t) ) + \uline\beta(t) \nabla u(t) + \gamma(t)\, u(t) 
            = 0 \mbox{ for } t\in (0,T],  
            \quad u(0)=u_0,
\end{equation}
where $u_0$ is the payoff and the coefficient functions are given by
\begin{equation}
    \uuline{\alpha}(t) := \frac{1}{2} \begin{pmatrix} \nu_{t} & \nu_{t} \sigma \rho \\ \nu_{t} \sigma \rho & \nu_{t} \sigma^2  \end{pmatrix},
    \qquad
    \uline\beta(t) := - \begin{pmatrix}  r - \frac{1}{2} \nu_{t} - \frac{1}{2} \sigma \rho \\ \kappa \theta - \kappa \nu_{t} - \frac{1}{2} \sigma^2 \end{pmatrix}, 
    \qquad
    \gamma(t) := r.
\end{equation}
Truncating $x\in\R^{2}$ to some bounded domain $\Omega \subset \mathbb{R}^2$ and posing appropriate (truncation) boundary conditions yields
\begin{subequations}\label{initboundHeston}
    \begin{align}
             \frac{d}{dt} u(t,x)+ \mathcal{A}(t)\, u(t,x) &= 0 \mbox{ in } (0,T]\times \Omega, \\
        u(t,x) &= 0 \mbox{ on }  [0,T]\times \partial \Omega , \\
        u(0,x) &= u_0(x) \mbox{ on } \Omega, \label{initboundHeston:c}
    \end{align}
\end{subequations}
where the differential operator is given as $\mathcal{A}(t) u := - \Div(\uuline\alpha(t) \nabla u ) + \uline\beta(t) \nabla u + \gamma(t) u$,  i.e., $\mathcal{A}(t) \in \mathcal{L}(V,V')$ with $V := H_0^1(\Omega)$. It is readily seen that $\mathcal{A}(t)$ is elliptic if and only if $\nu_{t}>0$ and $\rho \in (-1,1)$, which we will always assume.

It is well-known that the bilinear form associated to $\mathcal{A}(t)$, i.e., $\langle \mathcal{A}(t)\, u, v\rangle_{V'\times V} = a(t; u,v)$ reads
\begin{equation*}
    a(t; u,v) := \int_{\Omega} \{ \uuline{\alpha}(t) \nabla u(x)\cdot \nabla v(x) + \uline\beta(t)\cdot \nabla u(x)\, v(x) + \gamma(t)\, u(x)\, v(x)\} dx.
\end{equation*}
If $(\uuline\alpha)_{ij}(t)$, $\uline\beta_i(t)$, $\gamma(t) \in L_{\infty}(\Omega)$, $i,j = 1,2$, and $u_0\in L_2(\Omega)=:H$, then there exist constants $0< M_{a}, \alpha_{a}<\infty$ and $0\le\lambda_{a}<\infty$ such that for all $u,v\in V$ it holds 
\begin{subequations}\label{Eq:Energy}
    \begin{align}
    |a(t;u,v)| &\leq M_a \| u \|_{V} \| v\|_{V}, \label{Eq:Energy:1} \\ 
    a(t;u,u) + \lambda_a \| u\|^2_{H} &\ge \alpha_a \| u \|^2_{V},  \label{Eq:Energy:2}
\end{align}
\end{subequations}
if $\mathcal{A}(t)$ is elliptic, see, e.g.\ \cite[Thm.\ 2, p.\ 318]{Evans}. Obviously, \eqref{Eq:Energy:1} means boundedness and \eqref{Eq:Energy:2} is a G\aa{}rding inequality.

\subsubsection*{Parameters}
In order to calibrate the model, one has to determine the model parameters, i.e., $\mu_{1}:=(\rho, \kappa, \theta, \sigma)$ based upon market data. Moreover, in particular for pricing, we would like to determine prices for different payoff functions, i.e., we consider the initial value $\mu_{0}:=u_0$ as a parameter function.

\subsection{Reduced Basis Methods (RBMs)}
\label{Sec:1.2}
Now, we recall the main features of Reduced Basis Methods (RBMs) for paramaterized partial differential equations (PPDEs) that are needed here. For more details, we refer to the surveys \cite{HaasdonkLuminy,RP} and the references therein. We assume that the variational form of the PDE is given by
\begin{equation}\label{eq:PPDE}
    u(\mu)\in\X:\quad b(\mu; u(\mu), v) = f(\mu; v)\quad \forall v\in \Y,
\end{equation}
where $\X$, $\Y$ are (infinite-dimensional) Hilbert spaces, $b:\cD\times\X\times\Y\to\R$ a parametric bilinear form and $f:\cD\times\Y\to\R$ is a parametric functional, $f(\mu;\cdot)\in\Y'$. The parameter space is denoted by $\cD$. Well-posedess of \eqref{eq:PPDE} is always assumed, which in particular means that an inf-sup-condition
\begin{equation}\label{eq:infsup}
    \inf_{w\in\X} \sup_{v\in\Y} \frac{b(\mu; w,v)}{\| w\|_{\X}\, \| v\|_{\Y}}\ge \beta_{b} >0
\end{equation}
is satisfied for all $\mu\in\cD$. 

The next step is the assumption that a stable finite discretization $\X^\cN\subset\X$, $\Y^\cN\subset\Y$ in the sense
\begin{equation}\label{eq:infsup:cN}
    \inf_{w^\cN\in\X^\cN} \sup_{v^\cN\in\Y^\cN} \frac{b(\mu; w^\cN,v^\cN)}{\| w^\cN\|_{\X}\, \| v^\cN\|_{\Y}}\ge \beta^\cN >0
\end{equation}
is available that is sufficiently fine so that the discrete Petrov-Galerkin solution $u^\cN(\mu)$ is a sufficiently good approximation to $u(\mu)$. Hence, $u^\cN(\mu)$ is called \emph{detailed} or \emph{truth} solution. We assume that $u^\cN(\mu)$ can be computed with $\cN$ operations and that this number $\cN$ is too large to be acceptable for realtime or multi-query computations, so that a reduced model is necessary.

In order to derive this reduced model, the idea is to split the computations into an offline and an online phase. During the offline phase, one may afford to use the detailed model in order to compute \emph{snapshots} $u^i:=u^\cN(\mu^i)$, for $i=1, \ldots , N$ and $N\ll\cN$ is much smaller so that the reduced model is feasable online even if the complexity is of the order $\cO(N^{3})$ since the stiffness matrix is in general densely populated.

The reduced solution $u_N(\mu)$ is then determined as the (Petrov-)Galerkin solution in the space $\X_N:=\Span\{ u^i: i=1,\ldots ,N\}$ with a stable reduced space $\Y_N\subset \Y^\cN$ of dimension $N$:
\begin{equation}\label{eq:RBsol}
    u_N(\mu)\in \X_N:\quad b(\mu; u_N(\mu), v_N) = f(\mu; v_N)\quad \forall v_N\in \Y_N.
\end{equation}
Stability is to be understood in a uniform inf-sup-sense, i.e., 
\begin{equation}\label{eq:infsup:N}
    \inf_{w_N\in\X_N} \sup_{v_N\in\Y_N} \frac{b(\mu; w_N,v_N)}{\| w_N\|_{\X}\, \| v_N\|_{\Y}}\ge \beta_{\text{LB}} >0,
\end{equation}
with $\beta_{\text{LB}}$ independent of $N$ as $N\to\infty$ as well as of $\mu$.

The snapshots defining the reduced space are determined by the (offline-)selection of parameter samples $S_{N}:=\{ \mu^{i}: i=1,\ldots ,N\}$ and those samples are computed by maximizing a computable error estimate $\Delta_{N}(\mu)$ w.r.t.\ the parameter $\mu$. This can be done e.g.\ by nonlinear optimization or a greedy method w.r.t.\ a finite so-called \emph{training set} $\Xi^{\text{train}}\subset\cD$. Such an error estimate can e.g.\ be derived as follows.

\begin{prop}\label{Prop:DeltaN}
Let $\X_{N}\subset\X^{\cN}$ and $u_{N}(\mu)\in \X_{N}$ be the solution of \eqref{eq:RBsol}.
Defining the \emph{residual} by $r_N(\mu; v) := f(\mu; v) - b(\mu; u_N(\mu),v) = b(\mu; u^{\cN}(\mu) - u_N(\mu),v)$ for all $v\in \Y^{\cN}$, we obtain the following error estimate
\begin{equation}
    \|u^{\cN}(\mu) - u_N(\mu)\|_{\X^{\cN}} 
        \leq \frac{1}{\beta^{\cN}} \|r_N(\mu)\|_{(\Y^\cN)'} =: \Delta_N(\mu),
        \label{errest}
\end{equation}
where $\beta^{\cN}$ denotes the inf-sup-constant \eqref{eq:infsup:cN} of $b(\mu; \cdot,\cdot)$ on $\X^{\cN}\times\Y^{\cN}$ (possibly with discrete norms $\|\cdot\|_{\X^\cN}$, $\|\cdot\|_{\Y^\cN}$).
\hfill\qed
\end{prop}

The ultimate efficiency aim is to realize an online complexity that is independent of $\cN$. The key for that realization is the assumption that bilinear form and right-hand side functional are separable w.r.t.\ the parameter, i.e., 
\begin{equation}\label{eq:affDec}
    b(\mu; w,v)  = \sum_{q=1}^{Q_b} \vartheta_q^b(\mu)\, b_q(w,v),
    \qquad
    f(\mu; v) = \sum_{q=1}^{Q_f} \vartheta_q^f(\mu) f_{q}(v),
\end{equation}
with functions $ \vartheta_q^b$, $ \vartheta_q^f:\cD\to\R$ and parameter-independent forms $b_{q}:\X\times\Y\to\R$ and $f_{q}:\Y\to\R$. In the RB literature this is often --a bit misleading-- called \emph{affine decomposition}.

Let $\{ v^j: 1\le j\le N\}$ be a basis of the space $\Y_N$, then \eqref{eq:RBsol} reads
\begin{equation}\label{eq:RBsol-dis}
    \bu_N = 
    (u_{N,i})_{i=1,\ldots ,N}\in\R^N: \, 
        \sum_{i=1}^N u_{N,i}\, b(\mu; u^i, v^j) = f(\mu; v^j)\quad \forall j=1,\ldots ,N.
\end{equation}
Furthermore, let $u^{i}=\sum_{n=1}^{\cN} \alpha^{i}_n \varphi_{n}^\cN$ be the representation of the snapshots in a basis $\{\varphi_n^\cN:\, n=1,\ldots ,\cN\}$ of $\X^\cN$ and correspondingly $v^{j}=\sum_{n=1}^{\cN} \beta^{j}_n \psi_{n}^\cN$ in a basis $\{\psi_n^\cN:\, n=1,\ldots ,\cN\}$ of $\Y^\cN$. If \eqref{eq:affDec} holds, the computation of the stiffness matrix and the right-hand side of the reduced linear system can be split in an offline/online fashion as follows (here only for the stiffnes matrix)
\begin{align*}
    b(\mu; u^i, v^j)
        &= \sum_{n,n'=1}^{\cN} \alpha^{i}_n \beta^{j}_{n'}\, b(\mu; \varphi_{n}^\cN,  \psi_{n'}^\cN) 
        = \sum_{n,n'=1}^{\cN} \alpha^{i}_n \beta^{j}_{n'}\, \sum_{q=1}^{Q_b} \vartheta_q^b(\mu)\, b_q(\varphi_{n}^\cN,\psi_{n'}^\cN) \\
        &= \sum_{q=1}^{Q_b} \vartheta_q^b(\mu) 
            \sum_{n,n'=1}^{\cN} \alpha^{i}_n \beta^{j}_{n'}\, b_q(\varphi_{n}^\cN,\psi_{n'}^\cN) 
            =: \sum_{q=1}^{Q_b} \vartheta_q^b(\mu) (b_q)_{i,j},
\end{align*}
where the terms $(b_q)_{i,j}$ are $\mu$-independent and can thus be computed offline. Online, for a new parameter $\mu\in\cD$, the stiffness matrix of the RB-system is obtained by
$$
    (\bB_N(\mu))_{i,j} = b(\mu; u^i, v^j) = \sum_{q=1}^{Q_b} \vartheta_q^b(\mu) (\bB_q)_{i,j},\qquad
    i,j=1,\ldots N,
$$
which requires $\cO(Q_b\, N^2)$ operations independent of $\cN$. The same procedure can be done for the right-hand side and also for the error estimator $\Delta_N(\mu)$ since one can derive a separation like in \eqref{eq:affDec} also for the residual $\|r_N(\mu)\|_{(\Y^\cN)'}$.

Starting from an elliptic PPDE, there are several RB methods for parabolic problems using a usual time-stepping approach, e.g.\ \cite{GP,HO}. We do not follow this path here, in particular since we aim at using the initial value as a parameter function. Standard RB error estimates for time-stepping-based methods involve sums of residuals of each time step, so that an initial error is possibly heavily amplified. This is one of the reasons, why we consider a space-time variational formulation of parabolic PPDEs that will lead us to a problem of the form \eqref{eq:PPDE}, where $\X$, $\Y$ will be Bochner spaces, i.e. involve space and time.

\mysection{Space-Time Formulation of Parabolic PPDEs}
\label{Sec:3}

We now review the space-time variational formulation and its discretization.

\subsection{Parabolic PPDEs}
\label{Sec:3.1}

The Heston model problem yields a parameterized parabolic initial-boundary value problem of the following form: Let $I:=(0,T)$ be the (open) time interval and $V\hookrightarrow H\hookrightarrow V'$ be a Gelfand triple of Hilbert spaces (e.g.\ $V=H^1_0(\Omega)$, $H=L_2(\Omega)$ for a bounded domain $\Omega\subset\R^d$). The parameters are assumed to take the form $\mu=(\mu_0,\mu_1)\in\cD:=\cD_0\times\cD_1$, where $\cD_0\subset H$ is a set of possible initial values and $\cD_1\subset\R^P$ is a parameter space of finite (or infinite) dimension. Then, given some $g(\mu; t)\in V'$, $t\in I$ a.e., we look for $u(\mu;t)\in V$, $t\in I$ a.e., such that
\begin{subequations}
    \label{Eq:PIVP}
    \begin{align}
            \kern-5pt
        \langle\dot{u}(\mu;t), \phi\rangle_{V'\times V} 
        + a(\mu_1; u(\mu;t), \phi) &= \langle g(\mu_1;t), \phi\rangle_{V'\times V}\,\, \forall\, \phi\in V,\, t\in I\, \text{a.e.},
            \kern-10pt
            \label{Eq:PIVP:1}\\
        u(\mu; 0) &= \mu_0\,\,\text{in}\, H,
            \label{Eq:PIVP:2}
    \end{align}
\end{subequations}
where $a(\mu_1; \cdot, \cdot): V\times V\to\R$ is a bounded bilinear form. Note, that $a(\cdot; \cdot, \cdot)$ and $g(\cdot; \cdot)$ are assumed to depend only on $\mu_1$, \emph{not} on $\mu_0$. As described above, $\mu_0$ is the initial value parameter, whereas $\mu_1$ contains the parameters occurring within the coefficients of bilinear form and right-hand side, e.g.\ parameters to be calibrated. We can easily extend our findings to time-dependent bilinear forms $a(\mu; t, \cdot,\cdot)$ as well. 

As already motivated in Section \ref{Sec:1.2}, we assume the following separability w.r.t.\ the parameter $\mu_1\in\cD_1$ for $\phi, \psi\in V$   and $t\in I$
\begin{equation}
        a(\mu_1; \phi, \varphi) = \sum_{q=1}^{Q_a} \vartheta_q^a(\mu_1)\, a_q(\phi, \varphi), 
            \qquad
        g(\mu_1; t) = \sum_{q=1}^{Q_g} \vartheta_q^g(\mu_1)\, g_q(t),
        \label{eq:aff-a-g}
\end{equation}
where $\vartheta_q^a,  \vartheta_q^g: \cD_1\to\R$, $a_q(\cdot, \cdot): V\times V\to\R$ and $g_q(t)\in V'$ are given.

\subsection{Variational Form} 
For existence and uniqueness of a solution to \eqref{Eq:PIVP} we assume that there exist constants $M_a, \alpha_a > 0$ and $\lambda_a \in \mathbb{R}$ such that for all $\phi,\psi \in V$ and all $\mu_1\in\cD_1$
\begin{align}
    \label{Eq:bound}
     |a(\mu_1; \phi,\psi)| &\leq M_a \|\phi\|_V \|\psi\|_V                         &&\mbox{(boundedness), } \\
    \label{Eq:garding}
     a(\mu_1; \psi,\psi) + \lambda_a \|\psi\|^2_H &\geq \alpha_a \|\psi\|^2_V &&\mbox{(G\aa{}rding inequality)},
\end{align}
see \eqref{Eq:Energy}. 
Of course, in general the constants $M_a$, $\lambda_a$ and $\alpha_a$ depend on $\mu_1$ -- for simplicity we consider upper respectively lower bounds that are parameter-independent.
For the space-time variational form of \eqref{Eq:PIVP} we define $\Z:= L_2(I;V)$ and consider as in \cite{SS} the following trial space 
\begin{align}
    \X := \left\{ w \in \Z :\, \dot{w} \in \Z' \right\} = L_2(I;V) \cap H^1(I;V'),
\end{align}
with the norm $\|w\|^2_{\X} := \|w\|^2_{\Z} + \|\dot{w}\|^2_{\Z'} + \|w(T)\|^2_H$, $w\in\X$. The test space is $\Y := \Z \times H$ and for every $v = (z,h)$ in $\Y$ the norm is defined by $\|v\|^2_{\Y} := \|z\|^2_{\Z} + \|h\|^2_H$. For $w \in \X$ and $v = (z,h) \in \Y$ we define
\begin{align}
    b(\mu_1; w,v) 
        &:= \int_I \left\langle \dot{w}(t),z(t) \right\rangle_{V'\times V} dt + \int_I a(\mu_1; w(t),z(t)) dt + (w(0),h)_H \notag \\
        &=: b_1(\mu_1; w,z) + (w(0), h)_H \label{Eq:b1}
\end{align}
and the right-hand side as
\begin{equation}\label{Eq:g}
    f(\mu; v) := \int_I \left\langle g(\mu_1; t),z(t) \right\rangle_{V'\times V} dt + (\mu_0,h)_H =: g_1(\mu_1; z) + (\mu_0,h)_H.
\end{equation}
The space-time variational formulation of \eqref{Eq:PIVP} is of the form \eqref{eq:PPDE}, i.e., 
\begin{equation}
    \mbox{find } u(\mu) \in \X \mbox{ such that } b(\mu_1; u(\mu),v) = f(\mu; v) \ \ \forall v\in \Y.\label{variation}
\end{equation}
Due to the separability of $a(\mu_1; \cdot, \cdot)$ and $g(\mu_1; \cdot)$ w.r.t.\ the parameter $\mu_1$, we get a corresponding separation for the bilinear form as well ($v=(z,h)$):
\begin{align}
    b(\mu_1; w,v) 
    &= \int_I \langle \dot{w}(t), z(t)\rangle_{V'\times V} dt + \sum_{q=1}^{Q_a} \vartheta_q^a(\mu_1) \int_I a_q(w(t), z(t))\, dt + (w(0), h)_H \notag\\
    &=: \sum_{q=1}^{Q_b} \vartheta_q^b(\mu_1)\, b_q(w,v) \label{Eq:baff}
\end{align}
with $Q_b=Q_a+1$, $\vartheta_q^b = \vartheta_q^a$, $1\le q\le Q_a$, $\vartheta_{Q^b}^b\equiv 1$ as well as 
$b_q(w,v) = \int_I a_q(w(t), z(t))\, dt$, $1\le q\le Q_a$, and $b_{Q_b}(w,v) := (w(0), h)_H + \int_I \langle \dot{w}(t), z(t)\rangle_{V'\times V} dt$. The situation is slightly different for the right-hand side. Let $\Psi:=\{ \psi_m: \, m\in\N\}$ be a (Riesz-)basis of $H$, then $\mu_0 = \sum_{m\in\N} \mu_{0, m}\, \psi_m$, so that
$$
    f(\mu; v) = \sum_{q=1}^{Q_g} \vartheta_q^g(\mu_1) \int_I \langle g_q(t), z(t)\rangle_{V'\times V} dt 
        + \sum_{m=1}^\infty \mu_{0,m}\, (\psi_m, h)_H,
        \quad v=(z,h),
$$
which is a separation w.r.t.\ the parameters, but with infinitely many terms. If $\mu_0$ has some finite expansion (even in terms of a different set of functions or also obtained by an approximation), we would get an obvious separation with $Q_f=Q_g+L$, where $L$ denotes the number of terms in such a finite expansion of $\mu_0$. We will come back to this point later, in particular w.r.t.\ an efficient offline-online treatment of the right-hand side.

We define space-dependent quantities\footnotemark[1]
$$
    M_e:= \sup_{w \in \X \setminus \{0\}} \frac{\|w(0)\|_H}{\|w\|_{\X}},
    \qquad
    \varrho := \sup_{0 \neq \phi \in V} \frac{\|\phi\|_H}{\|\phi\|_V}
$$
\footnotetext[1]{Typically, one can bound $M_e\le\sqrt{3}$.}
and form-dependent ones (recall \eqref{eq:infsup})
$$
    \beta_a^* :=  \inf_{\mu_1\in\cD_1} \inf_{\phi\in V} \sup_{\psi\in V} \frac{ a(\mu_1; \psi, \phi)}{\|\phi\|_V\, \|\psi\|_V}, \qquad
    \beta_b :=  \inf_{\mu_1\in\cD_1} \inf_{w\in\X} \sup_{v\in \Y} \frac{ b(\mu_1; w, v)}{\| w\|_\X\, \| v\|_\Y},
$$
as well as lower bounds for the inf-sup-constant of the bilinear form $b$:
\begin{subequations}
    \label{Eq:beta}
        \begin{align}
        \label{Eq:beta:1}
        \beta^{\text{LB}}_\text{coer} 
            = \beta^{\text{LB}}_\text{coer} (\alpha_a, \lambda_a, M_a)
            &:= \frac{\min\{ \min\{1, M_a^{-2}\} (\alpha_a-\lambda_a\varrho^2), 1 \} }{\sqrt{2 \max\{1, (\beta_a^*)^{-2}\} + M_e^2}},\\
        \label{Eq:beta:2}
        \beta^{\text{LB}}_\text{time} (T)
            = \beta^{\text{LB}}_\text{time} (T, \alpha_a, \lambda_a, M_a)
            &:= \frac{\beta^{\text{LB}}_\text{coer}(\alpha_a, 0, M_a+\lambda_a\varrho^2)\, e^{-2\lambda_a T}}{\sqrt{\max\{2, 1+2 \lambda_a^2 \varrho^4\}}}.
        \end{align}
\end{subequations}

\begin{prop}[{\cite[Prop.\ 2.2, Cor.\ 2.7]{UP2}},{\cite[Thm.\ 5.1]{SS}}]\label{Prop:3.1}
    Let $a(\mu_1; \cdot,\cdot)$ satisfy \eqref{Eq:bound} and \eqref{Eq:garding}. Then, we obtain the inf-sup lower bound
    $\beta_b \ge 
        \beta_b^{\rm{LB}}:=\max\{ \beta^{\text{\emph{LB}}}_\text{\emph{coer}}, \beta^{\text{\emph{LB}}}_\text{\emph{time}} (T)\}$.\hfill$\qed$
\end{prop}

\begin{rem}
(a) Note that $\beta^{\text{\emph{LB}}}_\text{\rm{coer}}$ does not depend on time. However, this estimate is only meaningful (i.e., positive) if $\alpha_a-\lambda_a\varrho^2>0$ which means that $a(\mu_1; \cdot, \cdot)$ is coercive.
    
(b) In case of homogeneous initial conditions using 
    \begin{equation}\label{eq:X0}
        \W:=\{ w\in\X: w(0)=0\}, 
    \end{equation}    
        the above estimate holds for the form $b_1$ introduced in \eqref{Eq:b1} with $M_e=0$.
    
(c) As an example, let us consider the heat equation, i.e., $A(\mu_1)\equiv -\Delta$ (no parameter dependence), $V=H^1_0(\Omega)$, $H=L_2(\Omega)$. Hence, in this case, parameters only appear in the right-hand side, i.e.,  $a(\mu_1; \cdot, \cdot) \equiv a(\cdot, \cdot)$ as well as $b(\mu_1; \cdot, \cdot)\equiv b(\cdot, \cdot)$ and $f(\mu; v) = f(\mu_0; v) = \int_I \langle g(t), z(t)\rangle_{V'\times V} dt + (\mu_0, h)_H$, $v=(z,h)$. Then, we use $\|\phi\|_V^2=a(v,v)$, i.e., $M_a=1$, $\lambda_a=0$, $\alpha_a = \beta_a^*=1$. As in \cite[Cor.\ 2.5]{UP2}, we have $\beta_b\ge 1$.
\end{rem}

\begin{rem}\label{Rem:nonlin}
    At least for the numerical realization, we will also consider (quadratically) nonlinear problems, even though well-posedness will not be investigated here. In this case $a(\mu_{1}; \cdot,\cdot,\cdot): V^{3}\to\mathbb{R}$ is a trilinear form and in \eqref{Eq:PIVP:1} we would have the term $a(\mu_{1}; u(\mu;t), u(\mu;t),\phi)$ instead. The reason for this consideration is twofold: (1) Several financial models involve nonlinearities; (2) Our subsequent numerical approach using the space-time variational form particularly allows the treatment of polynomial nonlinearities, see \S \ref{Sec:nonlinear}  below.
\end{rem}

\subsection{Discretization}\label{Discretisation}

For a linear problem, one could reduce \eqref{Eq:PIVP}, in particular \eqref{Eq:PIVP:2} to a homogeneous initial condition. However, since we aim at considering the initial condition (also) as a parameter and also have nonlinear problems in mind, we keep the inhomogeneous initial condition, so that we need to modify what has been proposed in \cite{UP2}, see also \cite{AndreevDiss,AndreevHeat}.

We use (also) finite elements to construct finite dimensional subspaces $\X^{\cN} \subset \X$, $\Y^{\cN} \subset \Y$ and aim at determining an approximation $u^{\cN}(\mu) \in \X^{\cN}$ of the solution $u(\mu) \in \X$ of \eqref{variation}. The parameter $\cN$ will contain dimension parameters for time, space and initial value discretization spaces. Concerning notation, a (calligraphic) superscript will always denote a large (offline, `truth', detailed) dimension, whereas a reduced dimension will later be indicated with an (non-calligraphic) index. Moreover, spaces written with `blackboard bold' letters (\verb+\mathbb+) contain space-time functions, whereas `normal' letters are used for functions in space \emph{or} time.

In order to introduce the discretization, we basically follow \cite{SS,UP2} and  note that
\begin{equation}
    \X = H^1(I) \otimes V, \qquad
    \Y = \Z\times H := L_2(I;V) \times H = (L_2(I) \otimes V) \times H.
\end{equation}
Let $\mathcal{T}_{\Space}^{\cJ}$ be a triangulation of the underlying space $\Omega$. For discretizing the function space $V$ we consider the finite subspace $V^{\cJ}$ spanned by a nodal basis $\{ \phi_1, \ldots, \phi_{\cJ} \}$ with respect to the triangulation $\mathcal{T}_{\Space}^{\cJ}$. 

For the finite-dimensional temporal subspaces $E^\cK \subset H^1_{\{ 0\}}(I):=\{ \phi\in H^1(I):\, \phi(0)=0\}$ and $F^\cK \subset L_2(I)$ consider the discretization of the interval $I$ given by $\mathcal{T}_{\Time}^{\cK} := \{ t^k = k\, \Dt:\,  0 \leq k \le \cK, \Dt := \frac{T}{\cK} \}$. The trial space $E^\cK$ is spanned by the piecewise linear functions $\{\sigma^1, \ldots, \sigma^\cK \}$ w.r.t.\ $\mathcal{T}_{\Time}^{\cK}$. For every $1\le k \le \cK-1$, we choose $\sigma^k$ as the hat function with nodes $t^{k-1}$, $t^{k}$ and $t^{k+1}$ and the remaining ones are defined by $\sigma^\cK := \frac{t-t^{\cK-1}}{\Dt} \chi_{[t^{\cK-1},t^\cK]}$ and $\sigma^0 := \frac{t^1-t}{\Dt} \chi_{[0,t^1]}$ (which will be needed later). The test space $F^\cK$ is chosen as $\Span \{\tau^1, \ldots, \tau^\cK \}$ with respect to $\mathcal{T}_{\Time}^{\cK}$ where $\tau^k\equiv \chi_{I^k}$ is the characteristic function on $I^k := (t^{k-1},t^k]$.

Finally, the trial (and test) space $H^\cM\subset H$ for the initial condition  is denoted by $H^{\cM}:=\Span\{ \psi_1, \ldots ,\psi_\cM\}\subset V \subset  H$\footnote{The reason to impose $\psi_\ell\in V$ will also become clear a little later.}. We keep $H^\cM$ and in particular $\cM$ arbitrary here and will detail possible choices later. In particular, we also allow for the case $H^\cM=V^\cJ$ ($\cM=\cJ$) and discuss advantages and disadvantages of this choice.

With these preparations at hand, the discrete approximation subspaces of $\X$ and $\Y$ are defined as ($\langle\sigma^{0}\rangle:=\Span\{ \sigma^{0}\}$)
\begin{subequations}
    \begin{align}
        \X^{\cN} 
            &:= (\langle\sigma^0\rangle\otimes H^\cM) \oplus (E^\cK \otimes V^{\cJ}) 
                =: \Q^\cM \oplus\W^\cI, 
                \quad \cI=\cK\cdot \cJ,
                \label{Eq:Xeta}\\
        \Y^{\cN}  &:= (F^\cK \otimes V^\cJ) \times H^\cM 
                =: \Z^\cI\times H^\cM\footnotemark[2]. 
    \end{align}
\end{subequations}
\footnotetext[2]{For the test space, we could also replace $H^{\cM}$ by some other space, possibly also with different dimension. Here, we want to keep this issue simple.}%
Since $\dim(E^\cK) = \dim(F^\cK) = \cK$, $\dim(H^\cM)=\cM$ and $\dim(V^\cJ) = \cJ$, we have $\dim(\X^{\cN})=\cM + \cK \cJ=:\cN=\cI+\cM=\dim(\Y^{\cN})$. In case the discretized versions of $\X$ and $\Y$ would have different dimensions, a least squares method has to be used, \cite{AndreevDiss}. Finally, note that $\W^\cI\subset\W$, see \eqref{eq:X0}. 
We obtain a discrete variational formulation of \eqref{variation}: Find $u^{\cN}(\mu) \in \X^{\cN}$ such that 
\begin{equation}\label{discrete2}
    b(\mu_1; u^{\cN}(\mu),v^{\cN}) = f(\mu; v^{\cN}) \ \ \ \forall v^{\cN} \in \Y^{\cN},
    \qquad \mu=(\mu_0,\mu_1)\in\cD,
\end{equation}
which corresponds to a linear system that can be detailed as
$$
    b_1(\mu_1; u^{\cN}(\mu), z^\cI) + ( (u^\cN(\mu))(0), h^\cM)_H 
            = g(\mu_1; z^\cI) + (\mu_0, h^\cM)_H
$$
for all $v^{\cN}=(z^\cI,h^\cM)\in\Y^{\cN}=\Z^{\cI}\times H^{\cM}$. 

 \subsubsection{Time-stepping}
Similar to \cite{UP2}, it is not difficult to see that the above discretization is equivalent to a time-stepping scheme. In fact, we have the splitting $u^\cN=q^\cM+w^\cI\in \Q^\cM\oplus\W^\cI=\X^\cN$, in particular $u^\cN(0)=q^\cM(0)$ and $w^\cI(0)=0$. In terms of the respective bases, we get the representations
\begin{align}
    u^{\cN} &:= q^\cM+w^\cI 
        = \sum_{m=1}^\cM q_m (\sigma^0\otimes\psi_m) 
            + \sum_{k=1}^\cK \sum_{i = 1}^{\cJ} 
                w_i^k (\sigma^k \otimes \phi_i) \in \X^{\cN}, 
\end{align}
with the coefficient vectors $\bq_\cM := (q_m)_{m=1,\ldots , \cM}$ for the initial value as well as $\bw_\cI := (w^k_i)_{i=1,\ldots , \cJ,\, k=1,\ldots , \cK} =: (\bw_\cI^k)_{k=1, \ldots ,\cK}$ and (recall $\cN=\cM+\cI$)
\begin{align}
    v^{\cN} &= ( z^{\cI}, h^{\cM}) 
        = \Big(\sum_{\ell=1}^\cK \sum_{j = 1}^{\cJ} z_{j}^\ell (\tau^\ell \otimes \phi_j) , 
            \sum_{m = 1}^{\cM} h_{m} \psi_m\Big) \in \Y^{\cN},
        \label{Eq:vdM}
\end{align}
with the coefficient vectors  $\bz_\cI := (z_j^\ell)_{j=1,\ldots ,\cJ,\, \ell=1,\ldots, \cK}= (\bz_{\cI}^{\ell})_{\ell=1,\ldots,\cK}$ and $\bh_\cM:=(h_m)_{m=1,\ldots, \cM}$. With these notations, we obtain
\begin{align*}
b_1(\mu_1; u^{\cN}, z^{\cI} ) 
    &= \int_I \langle \dot{u}^{\cN}(t) , z^{\cI} (t) \rangle_{V' \times V}  
        + a(\mu_1; u^{\cN}(t), z^{\cI}(t))\,dt  \\
    &\kern-35pt= \sum_{m=1}^\cM \sum_{\ell=1}^\cK  \sum_{j=1}^{\cJ}  
        q_m z_{j}^\ell 
         \int_I  \langle  \dot{\sigma}^0(t) \psi_m, \tau^\ell(t) \phi_j \rangle_{V' \times V}  
        + a(\mu_1; \sigma^0(t) \psi_m, \tau^\ell(t)  \phi_j)\, dt  
        \\
    &\kern-25pt+ \sum_{k = 1}^\cK \sum_{\ell=1}^\cK \sum_{i,j=1}^{\cJ}  
        w_i^k z_{j}^\ell 
         \int_I  \langle  \dot{\sigma}^k(t) \phi_i, \tau^\ell(t) \phi_j \rangle_{V' \times V}  
        + a(\mu_1; \sigma^k(t) \phi_i, \tau^\ell(t)  \phi_j) \, dt  \\
    &\kern-35pt= \sum_{m=1}^\cM \sum_{\ell=1}^\cK  \sum_{j=1}^{\cJ}  
            q_m z_{j}^\ell 
        \Big\{ (\dot{\sigma}^0, \tau^\ell)_{L_2(I)} (\psi_m, \phi_j)_H 
        + (\sigma^0, \tau^\ell)_{L_2(I)} a(\mu_1; \psi_m, \phi_j)\Big\}  \\
    &\kern-25pt+ \sum_{k,\ell = 1}^\cK  \sum_{i,j=1}^{\cJ} 
        w_i^k z_{j}^\ell 
        \Big\{ (\dot{\sigma}^k, \tau^\ell)_{L_2(I)} (\phi_i, \phi_j)_H 
            + (\sigma^k, \tau^\ell)_{L_2(I)} a(\mu_1; \phi_i, \phi_j)\Big\}.
\end{align*}
Note, that $a(\mu_1; \psi_m, \phi_j)$ is well-defined since we have assumed that $\psi_m\in V$. For $k \geq 0$ and $\ell \geq 1$ we have $( \dot\sigma^k, \tau^\ell)_{L_2(I)} = \delta_{k,\ell} - \delta_{k+1,\ell}$ and $(\sigma^k, \tau^\ell)_{L_2(I)} = \frac{\Dt}{2} ( \delta_{k,\ell} + \delta_{k+1,\ell})$, in particular $( \dot\sigma^0, \tau^\ell)_{L_2(I)} = -\delta_{1,\ell}$ and $(\sigma^0, \tau^\ell)_{L_2(I)} = \frac{\Dt}{2} \delta_{1,\ell}$. 

Similar to \cite{UP2}, we set 
\begin{equation}\label{Eq:B}
    \bB^\cI(\mu_1) 
        := \bN_{\text{time}}^{\cK} \otimes \bM_{\text{space}}^{\cJ} 
        + \bM_{\text{time}}^{\cK}\otimes \bA_{\text{space}}^{\cJ}(\mu_1) 
        \in \R^{\cK\cJ\times \cK\cJ} = \R^{\cI\times\cI},
\end{equation}
where the temporal matrices read $\bN_{\text{time}}^{\cK} := ( ( \dot\sigma^k, \tau^\ell)_{L_2(I)}  )_{k, \ell=1, \ldots, \cK}$, $\bM_{\text{time}}^{\cK} := ((\sigma^k, \tau^\ell)_{L_2(I)} )_{k,\ell=1,\ldots , \cK}$ and the spatial ones are $\bM_{\Space}^\cJ := ((\phi_i, \phi_j)_H)_{i,j = 1, \ldots, \cJ}$ and $\bA_{\Space}^\cJ(\mu_1) = (a(\mu_1; \phi_i, \phi_j))_{i,j = 1, \ldots, \cJ}$. The matrix $\bB^\cI(\mu_1)$ was used in \cite[(2.14)]{UP2} to describe and analyze the discretization in the case of homogeneous initial conditions. Then, we obtain
\begin{align}
b_1(\mu_1; u^{\cN}, z^{\cI} ) 
    &= \sum_{k, \ell=1}^\cK \sum_{i,j=1}^{\cJ} 
        w_i^k z_{j}^\ell ( \bB^\cI(\mu_1))_{(k,i), (\ell,j)} 
        \notag \\
    &\kern-20pt
        + \sum_{m=1}^\cM \sum_{\ell=1}^\cK \sum_{j=1}^{\cJ} 
        q_m z_{j}^\ell 
        \Big\{ (\dot{\sigma}^0, \tau^\ell)_{L_2(I)} (\psi_m, \phi_j)_H 
            + (\sigma^0, \tau^\ell)_{L_2(I)} a(\mu_1; \psi_m, \phi_j)\Big\} 
        \notag \\
    &\kern-30pt
    = \bw_\cI^T \bB^\cI(\mu_1)\bz_\cI
         + \sum_{m=1}^\cM \sum_{j=1}^{\cJ} 
             q_m z_{j}^1 \Big\{ (- (\psi_m, \phi_j)_H 
            + \frac\Dt2 a(\mu_1; \psi_m, \phi_j)\Big\} 
        \notag\\
    &\kern-30pt
    = \bw_\cI^T \bB^\cI(\mu_1)\bz_\cI
        + \bq_\cM^T\Big(-\bM^{\cM,\cJ}_{\rm{i/s}} 
        + \frac\Dt2 \bA^{\cM,\cJ}_{\rm{i/s}}(\mu_1)\Big) \bz_\cI^1 
        \notag\\
    &\kern-30pt
        =: \bw_\cI^T \bB^\cI(\mu_1)\bz_\cI 
        +  \bq_\cM^T \bC^{\cM,\cJ}_{\rm{i/s}}(\mu_1) \bz_\cI^1,
    \label{Eq:wTBv}
\end{align}
where the involved matrices $\bM^{\cM,\cJ}_{\rm{i/s}}, \bA^{\cM,\cJ}_{\rm{i/s}}(\mu_1)\in\R^{\cM\times\cJ}$ 
are defined for $m=1,\ldots, \cM$, $j=1, \ldots , \cJ$ as
\begin{equation}\label{eq:Mats}
    (\bM^{\cM,\cJ}_{\rm{i/s}})_{m,j} :=  (\psi_m, \phi_j)_H,
    \quad
     (\bA^{\cM,\cJ}_{\rm{i/s}}(\mu_1))_{m,j} :=  a(\mu_1; \psi_m, \phi_j)_H.
\end{equation}

We split the coefficient vector for the unknown $\bu_\cN = (\bq_\cM, \bw_\cI)^T$ and set $\bu_\cN^0:=\bq_\cM$ as well as $\bu_\cN^k:=\bw_\cI^k$, $k=1,\ldots , \cK$ in order to formulate the time-stepping scheme. Then, we obtain for fixed $\ell\geq 1$ and any $j \in \{1, \ldots, n_h \}$ 
\begin{align}                
    b_1(\mu_1; u^{\cN}, \tau^\ell \otimes \phi_j)= \notag \\
        &\kern-95pt= 
            \begin{cases}
                [\bM_{\Space}^\cJ  \bu^1_\cN - (\bM^{\cM,\cJ}_{\rm{i/s}})^T\bu^0_\cN 
                    + \frac\Dt2 \big( \bA_\Space^\cJ(\mu_1)\bu^1_\cN 
                    + (\bA^{\cM,\cJ}_{\rm{i/s}}(\mu_1))^T\bu^0_\cN]_j
                    & \kern-7pt \mbox{if } \ell=1, \\
                \Dt \big[ \bM_{\Space}^\cJ \frac1\Dt (\textbf{u}_\cN^\ell 
                                - \textbf{u}_\cN^{\ell-1}) 
                    + \bA_{\Space}^\cJ(\mu_1) 
                        \frac{1}{2} (\textbf{u}_\cN^\ell + \textbf{u}_\cN^{\ell-1})\big]_j, 
                    & \kern-7pt \mbox{if } \ell>1
            \end{cases}        
            \notag\\        
        &\kern-60pt=: \Dt \big[ \frac1\Dt\,
                \bM^\ell_\cN (\textbf{u}_\cN^\ell - \textbf{u}_\cN^{\ell-1}) 
            + \bA_\cN^\ell(\mu_1)\, \textbf{u}_\cN^{\ell-1/2} \big]_j, 
            \label{Eq:MlAl}
\end{align}
where $\bu^{\ell-1/2}_\cN := \frac12 (\bu^\ell_\cN + \bu^{\ell-1}_\cN)$. 
On the right-hand side we use a trapezoidal approximation
\begin{align}
g(\mu_1; \tau^\ell \otimes \phi_j) 
    &= \int_I \langle g(\mu_1; t), \tau^\ell \otimes \phi_j (t,\cdot) \rangle_{V' \times V} dt 
    = \int_{I} \langle g(\mu_1; t), \tau^\ell (t) \phi_j \rangle_{V' \times V} dt \label{rhs} \\
    &\approx \frac{\Dt}{2} \langle g(\mu_1; t^{\ell-1}) + g(\mu_1; t^\ell) , \phi_j \rangle_{V' \times V}
        =: \Dt\, (\bg_\cN^{\ell-1/2}(\mu_1))_j. 
        \notag
\end{align}

\subsubsection{Initial value approximation}
Let us now discuss the approximation $\bu^0_\cN$ (given by $\bq_{\cM}$) of the coefficients of the initial value $u(0)=\mu_0$. Since $u^\cN(0)=q^\cM(0)$, $q^\cM\in\Q^\cM$, we have $u^\cN(0)=\sum_{m=1}^\cM q_m\, (\sigma^0 \otimes \psi_m)(0) = \sum_{m=1}^\cM q_m \, \psi_m\in H^{\cM}$, so that for $h^\cM=\sum_{m=1}^\cM h_m\psi_m\in H^\cM$, we get
\begin{align*}
    (u^{\cN}(0), h^{\cM})_H 
    &= \sum_{m, m'=1}^\cM q_{m'} h_{m}  (\psi_{m'}, \psi_m)_H 
    =  \bq_\cM^T {\bM}_\Init^{\cM} \bh_\cM,
\end{align*}
where $\bh_\cM:=(h_m)_{1\le m\le \cM}$ and ${\bM}_\Init^{\cM} := ((\psi_{m'}, \psi_m)_H)_{m', m=1, \ldots , \cM}$. The right-hand side of the discretization of \eqref{Eq:PIVP:2} for obtaining an approximation of the initial condition reads for the same $h^\cM\in H^\cM$
$$
    (\mu_0, h^\cM)_H = \sum_{m=1}^\cM h_m\, (\mu_0, \psi_m)_H,
$$
which is only computable if $\mu_0$ admits a finite expression. We refer to Remark \ref{Rem:finiteMu} below for possible choices. 

At this point we will only assume that a discretization of $\mathcal{D}_0$ is given by $\cD_{0}^\cL =  \Span \{ B_1,\ldots, B_\cL\}$, i.e., we consider $\mu^\cL_0 = \sum_{\ell = 1}^\cL \mu_0^{\ell} B_{\ell}$ and obtain 
$$
    (\mu^\cL_0, h^\cM)_H 
    = \sum_{\ell=1}^\cL \sum_{m=1}^\cM \mu_0^\ell \, h_m\, (B_\ell, \psi_m)_H
    = \boldsymbol{\mu}_{0, \cL}^T \bN^{\cL,\cM} \bh_\cM,    
$$
where $\bN^{\cL,\cM}:= ( (B_\ell, \psi_m)_H)_{\ell=1,\ldots ,\cL;\, m=1,\ldots ,\cM}$, $\bmu_{0,\cL} = (\mu_{0}^\ell)_{\ell=1,\ldots,\cL}$.  

\subsubsection{Crank-Nicolson scheme}
Putting everything together, we obtain the following Crank-Nicolson scheme for computing $\bu_\cN=(\bu_\cN^0,\ldots , \bu_\cN^\cK)$ (recall \eqref{Eq:MlAl}):
\begin{framed}\vspace*{-2mm}
\begin{subequations}\label{Eq:CN}
    \begin{align}
    (\bM_{\rm{init}}^{\cM})^T \textbf{u}^0_\cN &= (\bN^{\cL,\cM})^T \boldsymbol{\mu}_{0,\cL},
                \label{Eq:CN:1} \\
     \frac{1}{\Dt}  \bM^k_\cN (\textbf{u}^k_\cN - \textbf{u}_\cN^{k-1}) 
            + \bA_\cN^k(\mu_1)\, \bu^{k-1/2}_{\cN}
        &=  \bg_\cN^{k-1/2}(\mu_1),\,\,\, k=1, \ldots, \cK.
                \label{Eq:CN:2}
    \end{align}
\end{subequations}\vspace*{-4mm}
\end{framed}

Note, that $\bM_{\rm{init}}^{\cM}$ is regular, so that $\bu^0_\cN = \bu^0_\cN(\mu_0)$ is uniquely defined and the discrete problem \eqref{Eq:CN} is obviously well-posed.

\begin{rem}\label{Rem:specialDis}
    Let us now discuss some relevant special cases.
    
    (a) Let $H^\cM=V^\cJ$, $\cM=\cJ$. In this case, we get in \eqref{eq:Mats} that $\bM^{\cM,\cJ}_{\rm{i/s}}\equiv \bM_\Space^\cJ$ and $\bA^{\cM,\cJ}_{\rm{i/s}}\equiv \bA_\Space^\cJ$. This means that we do not need to distinguish the cases $\ell=1$ and $\ell>1$ in \eqref{Eq:MlAl}, i.e., we get a \emph{standard Crank-Nicolson scheme} with initial value $\bu^0_\cN$. 
        Moreover, $\bM^\cM_{\rm{init}}=\bM^\cJ_\Space$ (which is symmetric and positive definite (s.p.d.)) and $\bN^{\cL,\cM}= \bN^{\cL,\cJ} := ( (B_{\ell}, \phi_j)_H )_{1\le\ell\le \cL; 1\le j\le \cJ}$, i.e., we obtain $\bu^0_\cN = (\bM^\cJ_\Space)^{-1} (\bN^{\cL,\cJ})^T \boldsymbol{\mu}_{0,\cL}\in\R^{\cJ}$, which can be used as initial value for the Crank-Nicolson scheme \eqref{Eq:CN:2}.
        
        (b) If $\mu_0$ can be represented as (or approximated by) $\mu_0 = \mu_0^\cM = \sum_{m=1}^\cM \mu_{0}^m\psi_m\in H^\cM\ne V^{\cJ}$, $\bmu_{0,\cM}=(\mu_0^m)_{m=1,\ldots,\cM}$, then $\cL=\cM$, $\cD_{0}^{\cL}=H^{\cM}$, $\bN^{\cL,\cM}=\bM^{\cM}_{\rm{init}}$, which is s.p.d., so that $\bu^0_\cN = \bmu_{0,\cM}\in\R^\cM$. In this case, we need to modify the first step of the Crank-Nicolson scheme as we do not need to solve a linear system in \eqref{Eq:CN:1}.
\end{rem}

\begin{rem}\label{Rem:finiteMu}
    In view of Remark \ref{Rem:specialDis} (b) above, let us describe further scenarios for the approximation of $\mu_0\approx \mu_{0}^\cL\in \cD_0^{\cL} \subset H$ that we have in mind:

        (1) If $\Xi = \{ \xi_\ell:\, \ell\in\N\}$ is a Riesz basis for $H$, then $\mu_0$ has an expansion in that basis, i.e., $\mu_0=\sum_{\ell\in\N} \mu^{\ell}_{0} \xi_\ell$. Then, $\Xi_\cL:=\{ \xi_1,\ldots , \xi_\cL\}$, $\cD_0^\cL:=\Span(\Xi_\cL)$ may be selected as the `most significant' parts of the infinite expansion, e.g.\ by an adaptive approximation. The above approximation is then obtained using the corresponding expansion coefficients $\mu^{1}_{0}, \ldots , \mu^{\cL}_{0}$ or approximations of them (if they cannot be computed exactly).
        
        (2) Sometimes, the specific structure of possible initial values is known from the particular application. Then, it might be realistic (as in fact for some payoff functions in option pricing) that $\mu_0\in\Span\{ B_1, \ldots , B_\cL\}={\cD_{0}^\cL}\ne H^{\cM}$. In this situation we can directly write $\mu_0=\sum_{\ell=1}^\cL \beta_\ell B_\ell$ and $\bu_0^\cM$ is given as in \eqref{Eq:CN:1}.
        
        (3) As mentioned in Remark \ref{Rem:specialDis} (b), our approach particularly allows to choose $\cD_{0}^\cL$ as $H^\cM$ -- even though this results in an immediate smoothing of $\mu_0\in H$, since $H^\cM \subset V\subsetneq H$.
\end{rem}    

\subsubsection{Parameter Separation}
\label{Sec:PCN}

For later reference, we now detail the specific para\-meter-dependence of the discrete variational formulation. Obviously, $\bu^0_\cN$ can be computed by \eqref{Eq:CN:1} in dependency of $\mu_0$ (or its approximation $\mu_0^\cL$) in the sense that $\mu_0^\cL$ uniquely determines $\bu^0_{\cN}$ --- $\mu_1$ is not required. This can be formulated as follows. Recall from \eqref{Eq:Xeta} that $\X^{\cN}=\Q^\cM \oplus\W^\cI$, in particular
$$
    \W^\cI
        =\{ w_\cN\in\X^\cN:\, w_\cN(0)=0\} 
        = \Span\{ \sigma^k\otimes\phi_i:\, k=1,\ldots,\cK,\, i=1,\ldots, \cJ\}, 
$$
compare \eqref{eq:X0}. Then, \eqref{discrete2} can be divided as follows:
\begin{subequations}\label{Eq:Online}
\begin{align}
    \xi^{\cM}(\mu_0) \in H^\cM:
    && (\xi^{\cM}(\mu_0), h^\cM)_H 
        &= (\mu_0, h^\cM)_H 
        && \forall h^\cM\in H^\cM, 
        \label{Eq:Online:1}\\
    w^\cI(\mu) \in \W^{\cI}:
    && b_1(\mu_1; w^\cI(\mu), z^\cI) 
        &=\breve{f}(q^\cM(\mu_0), \mu_1; z^\cI)
        &&\forall z^\cI\in\Z^\cI,
        \label{Eq:Online:2}
\end{align}
\end{subequations}
with the extension of the initial value $q^\cM(\mu_0) := \sigma^0\otimes \xi^{\cM}(\mu_0)\in\Q^\cM$ and the modified right-hand side $\breve{f}(q^\cM(\mu_0), \mu_1; z^\cI):= g_1(\mu_1; z^\cI) - b_1(\mu_1; q^\cM(\mu_0), z^\cI)$.
In matrix-vector form as \eqref{Eq:wTBv} the second equation \eqref{Eq:Online:2} reads
\begin{align*}
    \bB^\cI(\mu_1)^T \bw_\cI(\mu) 
        &= \breve\bbf(\bu_\cN^0(\mu_0), \mu_1)\\
        &\kern-20pt
            := (g_1(\mu_1;\phi_i))_{i=1,\ldots, \cJ} 
                - (\bC^{\cM,\cJ}_{\rm{i/s}}(\mu_1)^T \bu_\cN^0(\mu_0), 0, \ldots ,0)^T. 
\end{align*}
The arising coefficient vectors define functions
\begin{align*}
    w^\cI(\mu) &= \sum_{k=1}^\cK \sum_{i=1}^{\cJ} w_i^k(\mu)\, (\sigma^k\otimes\phi_i)\in\W^{\cI}, &
    q^\cM(\mu_0) &= \sum_{m=1}^{\cM} q_m(\mu_0)\, (\sigma^0\otimes\psi_m)\in\Q^{\cM},
\end{align*}
so that $u^\cN(\mu):= q^\cM(\mu_0) + w^\cI(\mu)\in\X^\cN$ is the desired approximate solution. 

We stress the fact that \eqref{Eq:Online} can also be interpreted as a separation. In fact, \eqref{Eq:Online:1} determines an approximation of the initial value $\mu_0$ -- independent of $\mu_1$, whereas the evolution is determined in \eqref{Eq:Online:2} and -- as we have seen -- could be realized for example in terms of a Crank-Nicolson scheme.

\subsubsection{Nonlinear equations}
\label{Sec:nonlinear}

Note that the above mentioned separation is \emph{not} a consequence of the fact that a \emph{linear} parabolic problem allows one to reduce non-homogeneous initial conditions to homogeneous ones. In fact, if $a(\mu; \cdot,\cdot,\cdot)$ would be a trilinear form that induces a space-time trilinear form $b_{1}(\mu_{1}; \cdot,\cdot,\cdot)$, the analogue of \eqref{Eq:Online:2} would read
\begin{align*}
    & b_1(\mu_1; w^\cI(\mu), w^\cI(\mu), z^\cI) 
        + b_1(\mu_1; q^\cM(\mu_0), w^\cI(\mu), z^\cI) 
        + b_1(\mu_1; w^\cI(\mu), q^\cM(\mu_0), z^\cI) 
        = \\
    &=g(\mu_1; z^\cI) 
            - b_1(\mu_1; q^\cM(\mu_0), q^\cM(\mu_0), z^\cI),
\end{align*}
i.e., the quadratic term is supplemented by two linear terms since $q^\cM(\mu_0)$ is known at this stage. This also shows how to extend this approach to polynomial nonlinearities. The only difference is that the parameter induced by the initial value approximation is only in the right-hand side for a linear problem, but also appears within the coefficients of the PDE for polynomial nonlinearities.

\subsection{Stability of the discretization}
We have to show well-posedness of the discretized system, i.e., continuity, surjectivity and an inf-sup-condition for the linear operator induced by the bilinear form $b(\mu_1; \cdot, \cdot)$ on $\X^{\cN}\times\Y^{\cN}$, see \eqref{eq:infsup:cN}. Continuity and surjectivity are (more or less) readily seen.

For the case of homogeneous initial conditions, i.e. for $b_1(\mu_1;\cdot,\cdot)$, inf-sup-stability was investigated in \cite{UP2} under the condition that the bilinear form $a(\mu_{1}; \cdot,\cdot)$ satisfies a G\r{a}rding inequality \eqref{Eq:garding}, see also Proposition \ref{Prop:3.1}. However, the inf-sup-constant may deteriorate in the presence of strong advective terms, which is the case in the considered option pricing models. In that case one might require some sort of (known) stabilization.

In the general case, using the splitting \eqref{Eq:Online}, we could deduce the stability of the discrete problem \eqref{discrete2} from known results for both sub-problems (initial value approximation and Crank-Nicolson iteration). However, as Proposition \ref{Prop:DeltaN} shows, we need an explicit estimate for the discrete inf-sup-constant, which is \emph{not} easily obtained from the combination of the sub-problems.

For the special case presented in Remark \ref{Rem:specialDis} (a) (i.e., $H^{\cM}=V^{\cJ}$), the inf-sup stability was investigated in \cite{AndreevDiss,AndreevHeat}. In that case and for our chosen time discretization we only have to ensure a standard CFL condition (cf.\ \cite[Prop.\ 2]{AndreevHeat}) in order to obtain inf-sup-stability and also an estimate for $\beta^{\cN}$ in \eqref{eq:infsup:cN}. 

\mysection{A Reduced Basis Method (RBM) for Parameter Functions}
\label{Sec:4}
Now, we consider a Reduced Basis (RB) approximation for the Crank-Nicolson interpretation of the discrete space-time problem. Recall that the space-time variational formulation leads to a Petrov-Galerkin problem so that the reduced problem takes the form  \eqref{eq:RBsol}, where the bilinear form $b(\mu_1; \cdot, \cdot)$ only depends on $\mu_1$, whereas the right-hand side $f(\mu; \cdot)$ depends on the full parameter $\mu=(\mu_0, \mu_1)$.

As already pointed out, the parameter $\mu_0$ is a function. We are now going to describe a method to handle this challenge. For $\mu := (\mu_0, \mu_1) \in \cD_{0} \times \mathcal{D}_1$ the residual reads 
\begin{align*}
r_N(\mu;v)
        &= f(\mu; v) - b(\mu; u_N(\mu), v) \\
        &= g_1(\mu_1; z) + (\mu_0, h)_H - b_1(\mu_1; u_N(\mu), z) + (u_N(\mu)(0), h)_H \\
        &= g_1(\mu_1; z) - b_1(\mu_1; u_N(\mu), z) + (\mu_0-(u_N(\mu))(0), h)_H\\
        &=: r_{N,1}(\mu; z) + r_{N,0}(\mu; h),
\end{align*}
for any $v = (z, h) \in \Y^{\cN}$. Recall, that we need to construct a reduced basis that ensures a small residuum for the full parameter space $\cD$. In order to do so, we need an efficient online computation of the error estimator $\Delta_N(\mu)$ in \eqref{errest}, which requires a separation of the residual $r_N(\mu; v)$ into parts that depend only on $\mu$ and others depending only on $v$. This is no problem for $r_{N,1}(\mu; z)$ due to the separation properties of $g_1$ and $b_1$. However, the term $(\mu_0,h)_H$ is an issue since the inner product involves the parameter $\mu_0\in\cD_0$ and would be needed to be computed online e.g.\ in terms of a possibly costly quadrature.

\subsection{A two-step greedy method}
\label{Sec:2StepGreedy}
To construct a reduced basis, first assume that $\mu_0$ is, or can at least be well approximated by, a finite sum, i.e., $\mu_0 \approx \mu_0^\cL\in\cD_0^\cL=\Span\{ B_1,\ldots, B_\cL\}$. Then, we would get an usual separation of the form $(\mu_0,h)_H = \sum_{\ell=0}^\cL \beta_\ell(\mu_0)\, (B_\ell, h)_H$, where the terms $(B_\ell,h)_H$ can be precomputed offline. Constructing a basis with a standard greedy procedure would be possible by considering an $\cL$-dimensional parameter space containing the coefficients $\beta_\ell$. However, if $\cL$ is large, this is infeasable, also since the coefficients cannot easily be bounded, so that we would need to work with a $\cL$-dimensional hypercube of large `side lengths'. Consequently, determining parameter samples for snapshots e.g.\ by a greedy method might be extremely costly. Moreover, the use of $\cD_0^\cL$ might be a severe restriction to possible choices of the initial value $\mu_0$.

Hence, we need an alternative and consider again the residuum. The following estimate is immediate 
\begin{align*}
    \| r_N(\mu)\|_{\Y'}
    &= \sup_{v\in\Y} \frac{r_N(\mu; v)}{\| v\|_{\Y}} 
    = \sup_{(z, h)\in\Y^{\cN}} \frac{r_{N,1}(\mu; z) + r_{N,0}(\mu; h)}{( \| z \|_{\Z}^2 + \| h \|_H^2)^{1/2}} \\
    &= \| r_{N,1}(\mu)\|_{\Z'} + \sup_{h \in \mathcal{H}} \frac{r_{N,0}(\mu; h)}{\| h \|_H}
    \le \| r_{N,1}(\mu)\|_{\Z'} + \| \mu_0 - (u_N(\mu))(0)\|_H\\
    &=: R_{N,1}(\mu) + R_{N,0}(\mu_0).
\end{align*}
At a first glance it seems that $r_{N,1}$ (and $R_{N,1}$) only depends on $\mu_1$. However, the RB solution $u_N(\mu)$ involves both $\mu_0$ and $\mu_1$ so that both parameters enter. As already said earlier, the approximation of the initial value (and hence $R_{N,0}$), however, depends only on $\mu_0$.

The above form of the error estimate suggests to compute parameter samples (and snapshots) in a two-stage-method, namely first to determine samples $\mu_0^i$ for the initial value by maximizing $R_{N,0}(\mu_0)$ w.r.t.\ $\mu_0$ and second to consider the evolution and compute samples $\mu^{j}$ by maximizing $R_{N,1}(\mu^j)$ using the before-computed snapshots $h^\cM(\mu_0^i)$. This corresponds to the separated computations already introduced in \S \ref{Sec:PCN}. Let us now describe the two parts in detail.

We remark that even though we describe a greedy method, one could also use a different method to determine appropriate parameter samples e.g.\ by using nonlinear optimization w.r.t.\ the error estimate, \cite{MR2452388,UVZ}. The separation approach is independent of the particular maximization strategy.

\subsubsection*{Initial value greedy}
The first step is to generate a reduced basis for the initial value, i.e., we need the solution at $t=0$, which only depends on the parameter function $\mu_0 \in \mathcal{D}_0$, as we have seen in \S \ref{Sec:PCN}. For a given tolerance $\rm{tol}_0>0$, we are looking for $N_0$ samples $S^0_{N_0}:=\{ \mu_0^1, \ldots , \mu_0^{N_0}\}$ and corresponding snapshots $\{ h^1, \ldots, h^{N_0} \}$, such that $q^i := \sigma^0\otimes h^i\in\Q^\cM$ and $u_0^i:=q^i + 0 \in\X^\cN$ (i.e., $u_0^i:=q^i + w^i$, $w^i\in\W^\cI$, $w^i \equiv 0$, $u^i_0(0)=h^i$) is the corresponding snapshot.

Given a specific value of the parameter  $\mu_0\in \cD_0$, the corresponding (detailed) snapshot $h^\cM(\mu_0) \in H^\cM$ is determined by 
\begin{equation}\label{Eq:snap-init}
	(h^\cM(\mu_0), \psi_m)_H = (\mu_0,  \psi_m)_H,\qquad 1\le m\le\cM,
\end{equation}
where $\{ \psi_1,\ldots , \psi_\cM\}$ is the chosen basis for $H^\cM$. Set $u^i_0(0)=h^i:=h^\cM(\mu_0^i) = \sum_{n=1}^N \alpha_{0,n}(\mu_0) h^i$.

Given $h^1, \ldots ,h^N$ (where we should have $N\le N_0\ll\cM$) computed as snapshots corresponding to $S_N^0$, a corresponding RB initial value approximation $h_N(\mu_0)$ of $\mu_0\in\cD_0$ is determined by solving the linear equation system corresponding to
\begin{equation}\label{Eq:RB-init}
		(h_N(\mu_0), h^i)_H = (\mu_0,  h^i)_H,\qquad 1\le i\le N,
\end{equation}
provided that the inner products $(\mu_0,  h^i)_H$ can be computed online efficient (i.e., with complexity independent of $\cM$). Then, the error contribution reads
$$
	R_{N,0}(\mu_0) = \| \mu_0 - h_N(\mu_0)\|_H.
$$

One option to determine the reduced basis $h^1,\ldots ,h^{N_0}$ could be as follows: Choose $\eta^1,\ldots ,\eta^{\tilde{N}} \in \cD_0$, $\tilde{N}>N_0$, arbitrary, compute the Gramian $\bM^{\tilde{N}}_H := \big( (\eta^n, \eta^{n'})_H \big)_{n,n'=1,\ldots, \tilde{N}}$ and choose $\mu_0^1,\ldots , \mu_0^{N_0}$ as the (orthogonalized) eigenfunctions corresponding to the $N_0$ largest eigenvalues of $\bM^{\tilde{N}}_H$. This corresponds to a proper orthogonal decomposition (POD). If the $\eta^n$ are chosen well, this approach results in the best $H$-orthogonal choice. The obvious drawback is the strong dependence on the choice of the $\eta^n$.

For a greedy procedure, one chooses a training set $M^0_{\text{train}}\subset\mathcal{D}_0$ and determines parameter samples by maximizing $R_{N,0}(\mu_0)$ over $\mu_0\in M^0_{\text{train}}$. We obtain the greedy scheme in Algorithm \ref{Alg:GreedyInit}. The most crucial part may be to find a good training set $M_{\rm{train}}^0$. 

\begin{algo}
    \caption{Initial value greedy\label{Alg:GreedyInit}}
    \begin{algorithmic}[1]
    \STATE{Let $M_{\rm{train}}^0 \subset \mathcal{D}_0$ be the training set of initial values, $\rm{tol}_0>0$ a given tolerance.}
    \STATE{Choose $\mu_0^1 \in M^0_{\rm{train}}$, $S_1^0:= \{ \mu_0^1\}$, compute $h^\cM (\mu_0^1)$ as in \eqref{Eq:snap-init}, $\Xi_1^0 := \{ h^\cM(\mu_0^1) \}$ }
    \FOR{$j = 1, \ldots , N_0^{\rm{max}} $ }
        \STATE{$\mu_0^{j+1} = \arg\max\limits_{\mu_0\in M_{\rm{train}}^0} R_{j,0}(\mu_0) $}
        \STATE {\textbf{if} $R_{j,0}(\mu_0^{j+1})<\rm{tol}_0$ \textbf{then} $N_0:=j$, $H_{N_0}:=\Span(\Xi_{N_0}^0)$; \textbf{Stop} \textbf{end if}}
        \STATE{Compute $h^\cM(\mu_0^{j+1}) \in H^\cM$ as in \eqref{Eq:snap-init}.}
        \STATE{$S_{j+1}^0 := S_j^0 \cup \{ \mu_0^{j+1}\}$, $\Xi_{j+1}^0 := \Xi_j^0 \cup \{ h^\cM(\mu_0^{j+1}) \}$, orthogonalize $\Xi_{j+1}^0$.}
    \ENDFOR
    \end{algorithmic}
\end{algo}

It remains to discuss the efficient computation of the error term $R_{j,0}(\mu_0)$ for a given parameter $\mu_0\in\cD_0$. Note, that we obtain a set of orthonormal functions $\Xi_{N_0}^0$ as an output of Algorithm \ref{Alg:GreedyInit}. Hence, the RB approximation $h_{N_0}(\mu_0)$ coincides with the $H$-orthogonal projection of $\mu_0$ to $H_{N_0}$. This means that $R_{N_0,0}(\mu_0)$ is the error of the best approximation of $\mu_0$ in $H_{N_0}$. There are different possibilities to compute this error:
\begin{compactenum}
	\item If $\mu_0$ is given as formula, then an efficient quadrature may be used.
	\item If $\mu_0$ has a finite expression (like $\mu_0^\cL$ above) in terms of a stable basis $\{ B_1,\ldots ,B_\cL\}$, one may either use an efficient quadrature or transform $h_N(\mu_0)$ into that basis and use the coefficients of the difference.
	\item One could compute an orthonomal basis for the complement $H^\cM\ominus H_{N_0}$ and approximate $R_{N,0}(\mu_0)$ by computing coefficients of $\mu_0$ w.r.t.\ that complement basis (e.g.\ in terms of wavelets).
\end{compactenum}

\subsubsection*{Evolution greedy}
The next step is to find a basis for the part of the solution $u$ in $\W^\cI$, given the already determined reduced space $H_{N_0}$. Given a parameter $\mu=(\mu_0,\mu_1)\in\cD$ and an approximation $h^\cM(\mu_0)$, the evolution part $w^\cI(\mu) \in \W^{\cI}$ is computed as
\begin{equation}\label{Eq:snap-evol}
    b_1(\mu_1; w^\cI(\mu), z^\cI) 
        =g_1(\mu_1; z^\cI) - b_1(\mu_1; \sigma^0\otimes h^{\cM}(\mu_0), z^\cI)
        \quad\forall z^\cI\in\Z^\cI.
\end{equation}
For a reduced basis approximation corresponding to  $\mu=(\mu_0,\mu_1)\in\cD$, first compute 
$$
	h_{N_0}(\mu_0)=\sum_{n=1}^{N_0} \alpha_{0,n}(\mu_0)\, h^n
$$ 
as above. Then, given parameter samples $S_{N_1}^1=\{ \mu^1, \ldots , \mu^{N_1}\}\in\cD=\cD_0\times\cD_1$ (to be determined e.g.\ by a second greedy described below) and corresponding snapshots $w^{i} := w^\cI(\mu^i)\in\W^\cI$, a reduced basis approximation $w_{N_1}(\mu)\in \W_{N_1}=\Span\{ w^{i}: \, 1\le i\le N_1\}$ is determined by
\begin{equation}\label{Eq:RB-evol}
    b_1(\mu_1; w_{N_1}(\mu), z_{N_1}) 
    	= \breve{f}(\mu; z_{N_1})
       := g_1(\mu_1; z_{N_1}) - b_1(\mu_1; \sigma^0\otimes h_{N_0}(\mu_0), z_{N_1}),
\end{equation}
for all $z_{N_1}\in \Z_{N_1}$, 
where $\Z_{N_1}$ is a stable reduced space corresponding to $\W_{N_1}$ w.r.t.\ the inner product $b_1$ in the sense that
\begin{equation}\label{LBB:N}
	\inf_{w_{N_1}\in \W_{N_1}} \sup_{z_{N_1}\in \Z_{N_1}} \frac{b_1(\mu_1; w_{N_1}, z_{N_1})}{\| w_{N_1}\|_{\W}\, \| z_{N_1}\|_\Z} \ge \beta_1(\mu_1) >0 
\end{equation}
independent of $N_1\to\infty$, see \S \ref{sec:LBB} below. Here $\beta_1$ is the inf-sup constant of the bilinear form $b_1$.

It is readily seen that the right-hand side of \eqref{Eq:RB-evol} admits a separation w.r.t.\ the parameter. In fact, recalling \eqref{eq:aff-a-g} and \eqref{Eq:baff} (where here, as opposed to \eqref{eq:aff-a-g} we set  $b_{Q_b}(w,v) :=  \int_I \langle \dot{w}(t), z(t)\rangle_{V'\times V} dt$), we have
\begin{align}
	\breve{f}(\mu; z_{N_1})
	&:= g_1(\mu_1; z_{N_1}) - b_1(\mu_1; \sigma^0\otimes h_{N_0}(\mu_0), z_{N_1}) \notag \\
	&= \sum_{q=1}^{Q_g} \vartheta_q^g(\mu_1)\, g_q(z_{N_1})
		+ \sum_{q=1}^{Q_b} \vartheta_q^b(\mu_1)\, b_q(\sigma^0\otimes h_{N_0}(\mu_0),z_{N_1}) \notag \\
	&=  \sum_{q=1}^{Q_g} \vartheta_q^g(\mu_1)\, g_q(z_{N_1})
		+ \sum_{q=1}^{Q_b} \sum_{n=1}^{N_0} \vartheta_q^b(\mu_1)\, \alpha_{0,n}(\mu_0)\,
		b_q(\sigma^0\otimes h^n,z_{N_1}) \notag \\
	&=: \sum_{q=1}^{Q_g+N_0\, Q_b} \vartheta_q^{\breve{f}}(\mu)\, \breve{f}(z_{N_1})
		\label{eq:aff-fbreve}
\end{align}
with obvious definitions of the involved terms. Hence, we obtain an efficient offline-online splitting both for the computation of the reduced basis approximation $w_{N_1}(\mu)$ and of the residual $r_{N,1}(\mu;z)=g_1(\mu_1;z)-b_1(\mu_1;u_N(\mu),z)$, where we set $u_N(\mu):=q_{N_0}(\mu_0)+w_{N_1}(\mu)=\sigma^0\otimes h_{N_0}(\mu_0)+w_{N_1}(\mu)$, which means that
\begin{align*}
	r_{{N_1},1}(\mu;z)
	&= g_1(\mu_1;z)-b_1(\mu_1;u_N(\mu),z) \\
	&= g_1(\mu_1;z)-b_1(\mu_1;\sigma^0\otimes h_{N_0}(\mu_0),z) -b_1(\mu_1, w_{N_1}(\mu),z)\\
	&= \breve{f}(\mu; z)-b_1(\mu_1, w_{N_1}(\mu),z),
\end{align*}
which coincides with the residual of \eqref{Eq:RB-evol}. Recalling that $\cD_1\subset\R^P$, \eqref{Eq:RB-evol} is a reduced problem with a $(P+N_0)$-dimensional parameter space since $H_{N_0}$ is the RB initial value space. Such a dimension might be a challenge.  The error estimator is given by
\begin{equation*}
	\Delta_{N_1}^1(\mu) := \frac{\| r_{N,1}(\mu) \|_{\Z'}}{\beta_{\text{LB}}} = \frac{R_{N,1}(\mu)}{\beta_{\text{LB}}},
\end{equation*}
where $\beta_{\text{LB}}$ is a lower bound of the inf-sup constant of the bilinear form $b$, and we obtain a -- more or less -- standard greedy scheme described in Algorithm \ref{Alg:GreedyEvol}.

\begin{algo}
    \caption{Evolution greedy\label{Alg:GreedyEvol}}
    \begin{algorithmic}[1]
    \STATE{Let $M_{\rm{train}} \subset \cD$ be the training set, $\rm{tol}_1>0$ a given tolerance.}
    \STATE{Choose $\mu^1 \in M_{\rm{train}}$, $\mu^1:=(\mu_0^1, \mu_1^1)$, $S_1^1 := \{ \mu^1 \}$}
    \STATE{Compute the RB approximation $h_{N_0}(\mu_0^1) \in S^0_{N_0}$ as in \eqref{Eq:RB-init}}
    \STATE{Compute $w^\cI(\mu^1) \in \W^\cI$ as in \eqref{Eq:snap-evol}, $\Xi_1^1 = \{ w^\cI(\mu^1) \}$}
        \FOR{$j=1, \ldots , N_1^{\rm{max}} $ }
						\STATE{$\mu^{j+1} = \arg \max\limits_{\mu \in M_{\rm{train}}} \Delta_{j}^1( \mu)$}
            \STATE {\textbf{if} $\Delta_{j}^1(\mu^{j+1})<\rm{tol}_1$ \textbf{then} $N_1:=j$, $\W_{N_1} := \Span (\Xi_{N_1}^1)$; \textbf{Stop} \textbf{end if}}
    	    \STATE{Compute the RB approximation $h_{N_0}(\mu_0^{j+1}) \in S^0_{N_0}$ as in \eqref{Eq:RB-init}}
    	    \STATE{Compute $w^\cI(\mu^{j+1}) \in \W^\cI$ as in \eqref{Eq:snap-evol}}
            \STATE{$S^1_{j+1} := S_{j}^1 \cup \{ \mu^{j+1} \}$, $\Xi_{j+1}^1 := \Xi_{j}^1 \cup \{ w^\cI(\mu^{j+1})\}$}
        \ENDFOR
    \end{algorithmic}
\end{algo}

An obvious question arises how to choose the training set $M_{\rm{train}}$ in Algorithm \ref{Alg:GreedyEvol}, in particular the training samples for the initial value parameter. Possible choices for the the subset of $\cD_0$ include $M_{\rm{train}}^0$ from Algorithm \ref{Alg:GreedyInit} or --much smaller-- $S_{N_0}^0$, which might be a reasonable choice after performing already a greedy search for the initial value. We will come back to this point in our numerical experiments in Section \ref{Sec:5}.

\subsection{Stable RB test spaces}\label{sec:LBB}
It remains to construct a stable test space $\Z_{N_1}$ in the sense of \eqref{LBB:N}. It is well-known that it might be beneficial to construct this space by so-called \emph{supremizers} in an efficient offline-online manner, \cite{DPW,GV,RV}.

Let $\{ w^1,\ldots , w^{N_1}\}$ be the basis of $\W_{N_1}$ and fix $\mu_1 \in \mathcal{D}_1$. Then, the supremizer $s^n(\mu_1)\in\Z^\cI$, $1\le n\le N_1$, is defined by the relation
\begin{equation*}
    s^n(\mu_1) 
    	:= \arg\sup_{z^\cI \in \Z^{\cI}} \frac{b_1(\mu_1; w^n ,z^\cI)}{\|z^\cI\|_{\Z}}. 
\end{equation*}
In order to compute this quantity, recall that $\{ \zeta_i:= \tau^k\otimes\phi_j: 1\le k\le\cK, 1\le j\le\cJ, i=(k,j)\}$ is the basis of $\Z^\cI$, $\cI=\cK\cdot\cJ$, and note that 
\begin{align*}
    \| z^\cI\|_{\Z}^2 
    	&= \| z^\cI\|_{L_2(I;V)}^2 
	= \sum_{k,k'=1}^\cK \sum_{j,j'=1}^{\cJ} z_j^k z_{j'}^{k'}\, (\tau^k, \tau^{k'})_{L_2(I)}\, (\phi_j, \phi_{j'})_V  \\
    	&= \boldsymbol{z}^T_\cI (\boldsymbol{I}_{\Time}^{\cK} \otimes\boldsymbol{G}_{\Space}^\cJ) \boldsymbol{z}_\cI 
	=:  \boldsymbol{z}^T_\cI \boldsymbol{Z}^\cI \boldsymbol{z}_\cI
\end{align*}
with the Gramian matrices $\boldsymbol{G}_{\Space}^\cJ = ((\phi_i, \phi_j)_V)_{i,j=1,\ldots, \cJ}$ for $V^\cJ$ (w.r.t.\ the $V$-inner product) and $\boldsymbol{I}_{\Time}^{\cK} = ((\tau^k,\tau^{k'})_{k,k' = 1,\ldots,\cK} = (\Delta t)\, \boldsymbol{Id} \in\R^{\cK\times\cK}$. Next, let the expansion of $w^n$ in terms of the full basis $\{ \varpi_i:= \sigma^k\otimes\phi_j: i=(k,j), 1\le k\le\cK, 1\le j\le\cJ\}$ of $\W^\cI$ be denoted by
$$
	w^n = \sum_{i=1}^\cI \omega^n_i \varpi_i,
	\qquad
	\boldsymbol{\omega}^n_\cI := (\omega^n_i)_{i=1,\ldots,\cI}.
$$
Then, setting $z_i:= z_j^k$, $i=(k,j)$, we get
\begin{align*}
	b_1(\mu_1; w^n, z^\cI)
		&= \sum_{i,i'=1}^\cI \omega_i^n z_{i'} b_1(\mu_1; \varpi_i, \zeta_{i'})
		= (\boldsymbol{\omega}^n_\cI)^T \bB^\cI(\mu_1) \bz_\cI
		= (\bz_\cI)^T (\bB^\cI(\mu_1))^T \boldsymbol{\omega}^n_\cI.
\end{align*}
The vector $\bs^n(\mu_1)$ containing the expansion coefficients of $s^n(\mu_1)$ is given by
$$
	\bs^n(\mu_1) = ( \boldsymbol{Z}^\cI)^{-1} (\bB^\cI(\mu_1))^T \bw_\cI^n,
$$
In view of the separation of $b_1(\mu_1; \cdot, \cdot)$ w.r.t.\ the parameter $\mu_1$ (see \eqref{Eq:baff} with $b_{Q_b}$ replaced by $b_{Q_b}(w,v) := \int_I \langle \dot{w}(t), z(t)\rangle_{V'\times V} dt$), we have
$
	\boldsymbol{B}^{\cI}(\mu_1) 
        = \sum_{q=1}^{Q_b} \vartheta_q^b(\mu_1)\, \boldsymbol{B}_q^\cI
$ 
with parameter-independent matrices $\boldsymbol{B}_q^\cI$.  Hence, we obtain the representation 
$
\bs^n(\mu_1) = \sum_{q=1}^{Q_b} \vartheta_q^b(\mu_1) ( \boldsymbol{Z}^\cI)^{-1} (\bB^\cI_q)^T \bw_\cI^n
$
and the terms $\bz_q^n:= ( \boldsymbol{Z}^\cI)^{-1} (\bB^\cI_q)^T \bw_\cI^n$ can be computed offline (as they are parameter-independent). Since the $\mu_1$-dependent supremizers can be build by linear combinations (with $\mu_1$-dependent coefficients) of the functions $z_q^n\in\Z^\cI$ corresponding to the coefficient vectors  $\bz_q^n$, $1\le n\le N_1$, $1 \le q \le Q_b$ we choose for every $\mu_1 \in \cD_1$
$$
	\Z_{N_1}(\mu_1) := \Span \{ s^1(\mu_1), \ldots , s^{N_1}(\mu_1) \}
$$ 
as reduced test space, where $s^n(\mu_1) = \sum_{q = 1}^{Q_b} \vartheta_q^b(\mu_1)\, z^n_q$ .

\mysection{Numerical results}
\label{Sec:5}

\subsection{Heston Model} 
We consider the Heston model as described in \S \ref{Sec:1.1} above. The initial value $\mu_0$ corresponds to the payoff function $u_0$, see \eqref{initboundHeston:c}. For pricing problems, one aims at rapidly changing the payoff, which is the motivation to use the parameter function.

 In order to ensure well-posedness of the PDE, we require the natural assumptions $\nu_t>0$ for the volatility and $\rho\in (0,1)$ for  the correlation. Since we do not transform the initial conditions to homogeneous ones but work in  Bochner spaces using the space-time variational approach, we do not need additional conditions for the parameter spaces as e.g.\ in \cite{Hilber}.

Just for the ease of implementation, we choose homogeneous Dirichlet conditions. Unfortunately one has to work with a large domain $\Omega$ and with a fine discretization to get good results for the Crank Nicolson solution in comparison to the closed form solution of the Heston model (cf.\ \cite{Heston}) that one can use for validation. One could further improve the results by using e.g.\ the boundary conditions proposed in \cite{WAW}. 
\subsection{Initial condition parameter function} 
As already said earlier, we want to use the initial condition as a parameter function. Since payoff functions are not completely arbitrary, but have certain shapes, we introduce a model using Bernstein polynomials that allows for a small parametric representation of the payoff functions, which are usually continuous, piecewise smooth, convex and are composed of linear functions. 

\subsubsection{Bernstein Polynomials}
Bernstein polynomials are $H^1$-functions, preserve convexity and can be adapted locally, see, e.g.\  \cite{Sauer}.  We briefly recall the main properties that will be needed here. Consider an interval $\Delta = [v_0,v_1]\subset\R$, where $v_0< v_1$. Any $x\in\Delta$ has the unique representation 
\begin{equation*}
    x = u_0(x|\Delta)\, v_0 + u_1(x|\Delta)\, v_1,
\end{equation*}
where $u_j( \cdot |\Delta)\in\mathbb{P}_1$ (a linear polynomial) is nonnegative on $\Delta$ and $u_0(x |\Delta) +u_1(x |\Delta) = 1$ (convex combination, partition of unity). Denoting by
\begin{equation*}
    \Gamma_n := \big\{ \alpha = (\alpha_0, \alpha_1)\in\Z^2: \, \alpha_0, \alpha_1\ge 0, \alpha_0+\alpha_1= n\big\}
\end{equation*}
 the set of all homogeneous multi-indices of length $n$, the $\alpha$-th \emph{Bernstein-B\'{e}zier basis polynomial} of degree $n$ is defined for $\alpha\in\Gamma_n$ as
\begin{equation*}
    B_{\alpha}(x |\Delta) :=  \begin{pmatrix} n \\ \alpha \end{pmatrix} u_0^{\alpha_0}(x|\Delta) u_1^{\alpha_1} (x|\Delta).
\end{equation*}
The \emph{B\'{e}zier surface} for a given function $f: \Delta \to \mathbb{R}$ is defined by
\begin{equation*}
    B_{n}f(x|\Delta) = \sum_{\alpha \in \Gamma_n} f(x_{\alpha}) B_{\alpha} (x|\Delta),
    \qquad
    \text{where } x_{\alpha}\in\Delta \text{ s.t.\ } u(x_{\alpha} |\Delta) = \frac{\alpha}{n}.
\end{equation*}

\subsubsection{Payoff approximation for the Heston model}

For the Heston model, we consider a domain $\Omega = \Omega_1 \times \Omega_2 \subset \mathbb{R}^2$, where $\Omega_1$ models the logarithmic asset price $\log(S)$ and $\Omega_2$ the volatility $\nu$. The payoff depends only on $S$. Divide the interval $\exp(\Omega_1)$ into subintervals $I_1, \ldots, I_{\cL-1}$, $I_\ell=[v_\ell^-, v_\ell^+]$, $\ell= 1, \ldots, \cL-1$,  $v_{\ell-1}^+ = v_\ell^-$ for $\ell = 2,\ldots,\cL-1$, and choose $\Delta$ as each of these $I_\ell$.

We now model the initial value $\mu_0$ on each $I_\ell$ as a B\'ezier surface of degree $1$, i.e.,
\begin{align*}
    B_{1} \mu_0(x|I_\ell) 
    	&= \mu_0(x_{(1,0)}) B_{(1,0)}(x|I_\ell) + \mu_0(x_{(0,1)}) B_{(0,1)}(x|I_\ell) \\
	&= \mu_0(v_\ell^-)\, u_0^1 (x|I_\ell)\, u_1^0(x|I_\ell) + \mu_0(v_\ell^+)\, u_0^0 (x|I_\ell)\, u_1^1(x|I_\ell)\\
	&= \mu_0(v_\ell^-)\, u_0 (x|I_\ell) + \mu_0(v_\ell^+)\, u_1(x|I_\ell).
\end{align*}
The approximation of $\mu_0$ is then defined by
\begin{align*}
    \mu_0^\cL(x)  
    	&:= \sum_{\ell = 1}^{\cL-1} \{\mu_0({v_\ell^-}) u_{0}(x|I_\ell)  + \mu_0({v_\ell^+}) u_{1}(x|I_\ell)\}\\
   	&= \mu_0(v_1) u_{0}(x|I_1) 
    		+ \sum_{\ell = 2}^{\cL-1} \mu_0(v_\ell) (u_1(x|I_{\ell-1}) + u_0 (x|I_\ell))
		+ \mu_0(v_\cL) u_1(x|I_{\cL-1}),
\end{align*}
where we have used the fact that $\mu_0({v_{\ell-1}^+}) = \mu_0({v_\ell^-})$ and renamed ${v_\ell^-}$ by $v_\ell$ as well as $v_{\cL} := {v_{\cL-1}^+}$. Since $u_0({v_\ell}|I_\ell) = u_1({v_{\ell+1}}|I_\ell ) = 1$ and $u_0({v_{\ell+1}}|I_\ell) = u_1({v_\ell}|I_\ell ) = 0$ the following equality is obvious:
\begin{equation*}
    u_1(x|I_{\ell-1}) + u_0(x|I_\ell) 
    	= \frac{x - v_{\ell-1}}{v_\ell - v_{\ell-1}} \chi_{I_{\ell-1}}(x) 
		+ \frac{x - v_{\ell+1}}{v_\ell-v_{\ell+1}} \chi_{I_\ell}(x)=: B_\ell(x),
\end{equation*}
so that $\mu_0^\cL(x) = \sum_{\ell = 1}^{\cL} \mu_0(v_\ell) B_\ell(x)$.

\subsection{Numerical results}
We now present our numerical results.

\subsubsection{Data}
We used the following data for our simulation:
\begin{compactitem}
	\item We use the correlation as calibration parameter, i.e., $\mu_1 := \rho \in (0,1)=:\cD_1$;
	\item $\Omega_1 := \log([10^{-8},200])$ for the asset price, $\Omega_2 := (0,1]$ for the volatility, detailed dimension $\cJ = 14,271$;
	\item $\kappa = 0.8, \sigma = 0.6, \theta = 0.2$ and $r\equiv 0.001$;
	\item $T = 0.25$ (3 months), $\cK = 25$;
	\item $H^\cM := V^\cJ$ and $\cN = \cK\cdot\cJ + \cJ = \cJ(\cK+1) = 371,046$;
	\item Knots $\{v_1,\ldots,v_8\} = \ln(\{0,70,80,90,100,110,200\})$for the B\'{e}zier polynomials, 
			which implies that payoff functions with specific strike price $K \in \{70,80,90,100,110\}$ are represented exactly.
\end{compactitem}
All experiments have been performed on an Intel Core i7-3770 (3.40 GHz), with 16 GB RAM and using Matlab ($8.0.0.783$ (R2012b)). All RB calculations were implemented in RBmatlab, see \url{http://www.morepas.org}.

\subsubsection{Initial value reduced approximation}
\begin{figure}[!htb]
			\includegraphics[width=0.45\textwidth]{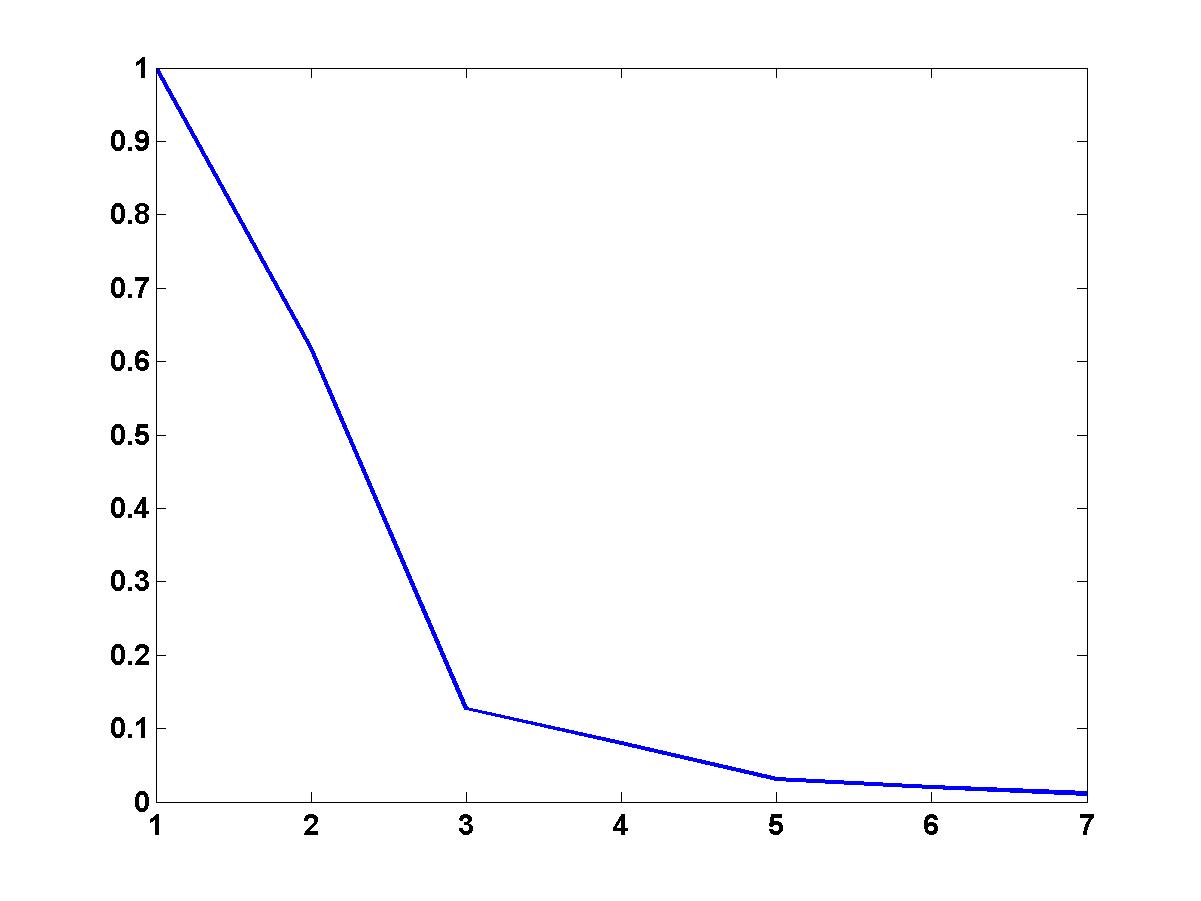}
			\caption{Eigenvalue decay for the Gramian matrix}
			\label{fig:EWdecay}
\end{figure}

\begin{figure}[!htb]
	\subfigure[First basis function $h^1$.]{
			\includegraphics[width=0.47\textwidth]{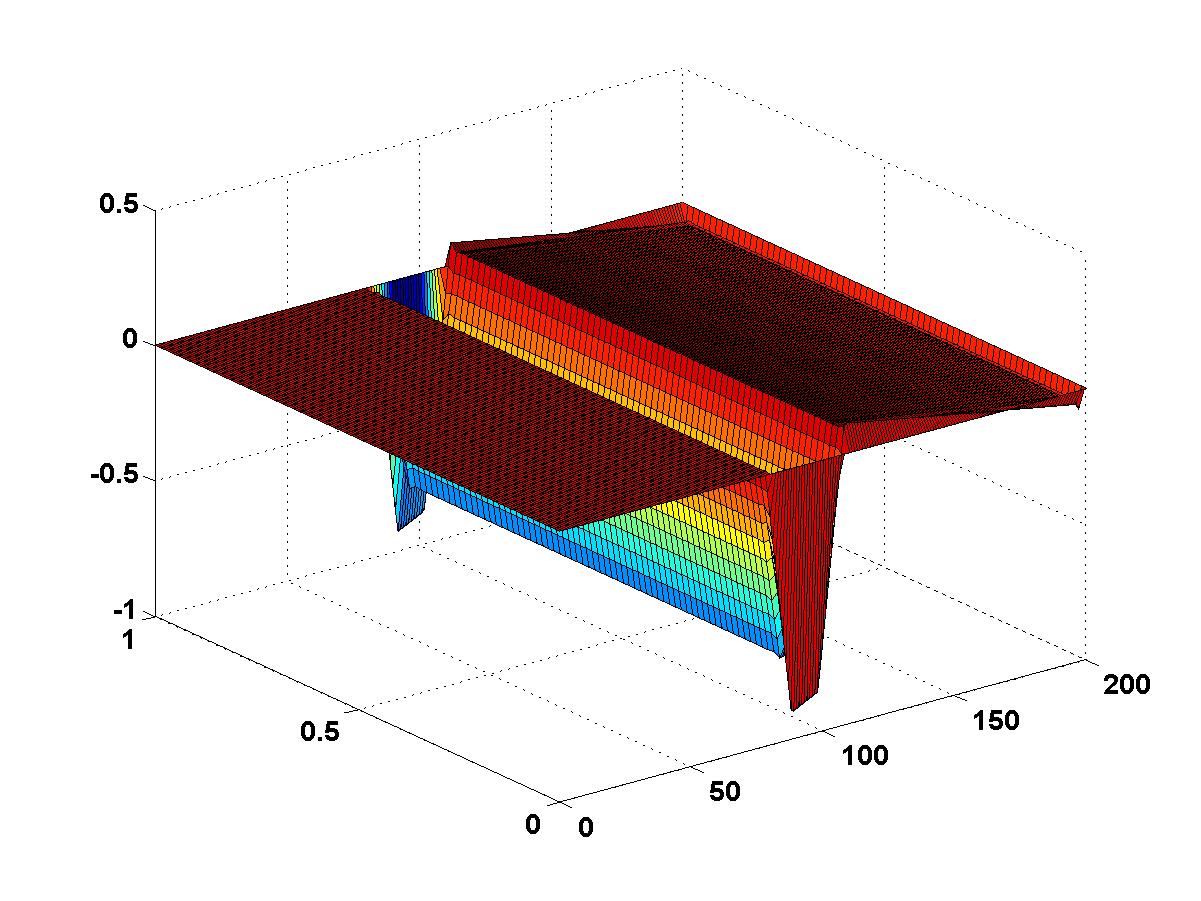}
			\label{fig:h1}}
	\subfigure[Second basis function $h^2$.]{
			\includegraphics[width=0.47\textwidth]{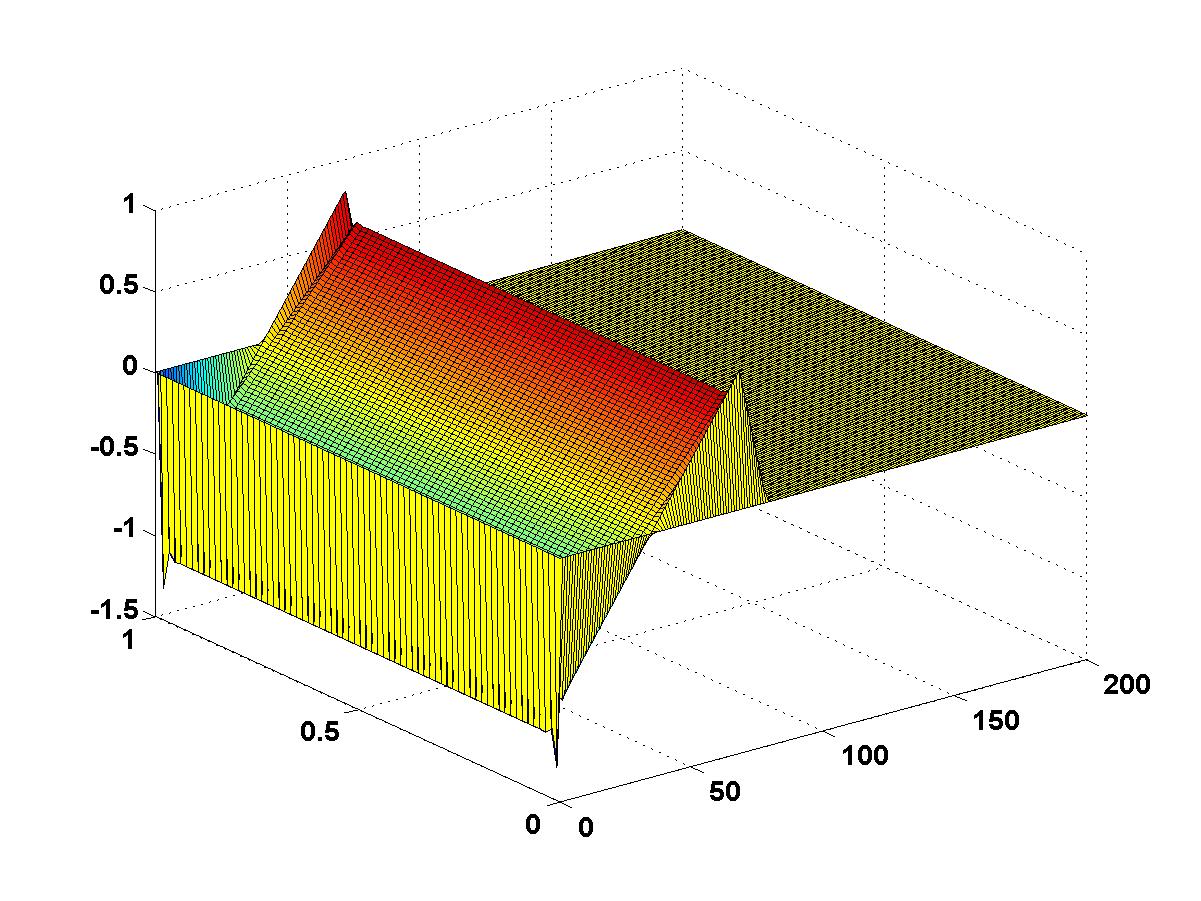}
			\label{fig:h2}}
	\caption{\label{fig:h}First two POD eigenfunctions projected to $V^\cJ$.}
\end{figure}

For determining the RB approximation of the initial value, we used a POD method based upon the linear Bernstein polynomials (in the notation of \S \ref{Sec:2StepGreedy} $\eta^\ell := B_{\ell}$, $\ell=1,\ldots , 7=:\tilde{N}$). We show the decay of the eigenvalues of the Gramian $\bM_H^{\tilde{N}}$ in Figure \ref{fig:EWdecay}. Choosing five basis functions $h^1,\ldots,h^5$ (i.e., a RB space $H_5$) results in a relative error of $0.0314$. We will investigate later how the choice of only $5$ POD basis functions influences the reduced solution of the whole problem. The first two orthogonal POD functions are shown in Figure \ref{fig:h}. 

For functions that can be modeled by a small number of basis functions as it is the case for call and put options using B\'{e}zier polynomials, it might be advisable to skip this first step and directly enter the evolution greedy using the small basis as a part of the training set.

\subsubsection{Evolution greedy}
In the second step we perform the evolution greedy in Algorithm \ref{Alg:GreedyEvol} to compute the RB space $\W_{N_1}$, where we first use $M_{\text{train}} :=  \{h^1,\ldots,h^5\} \times\{x_k = -0.5 + k\, \Delta s: k = 0,\ldots,11, \Delta s = \frac{1}{11} \} \subset \Span \{ B_1,\ldots,B_7\}  \times (-1,1)$ as a training set, i.e., $\# M_{\text{train}}=60$. 
As in \cite{UP2} we use a natural discrete space-time norm given by
\begin{equation*}
	\| w \|_{\bar{\X}}^2 := \| \bar{w} \|_{L_2(I;V)}^2 + \| \dot w\|_{L_2(I;V')}^2 + \| w(T)\|_H^2, \ w\in \X^\cN,
\end{equation*}
where $\bar{w}^k := (\Delta t)^{-1}\int_{I^k} w(t)\, dt\in V$ and $\bar{w}:=\sum_{k=1}^\cK \chi_{I^k}\otimes\bar w^k\in L_2(I; V)$. Using this norm, the discrete problem appearing in the evolution greedy is well-posed. Another reason for this choice is that the evolution greedy uses the same bilinear form as was used in \cite{UP2} to treat homogeneous initial conditions.

\begin{figure}[!htb]
	\subfigure[Maximum error over iterations of the strong evolution greedy and using the error estimate $R^1_{N_1}$]{
				\includegraphics[width=0.45\textwidth]{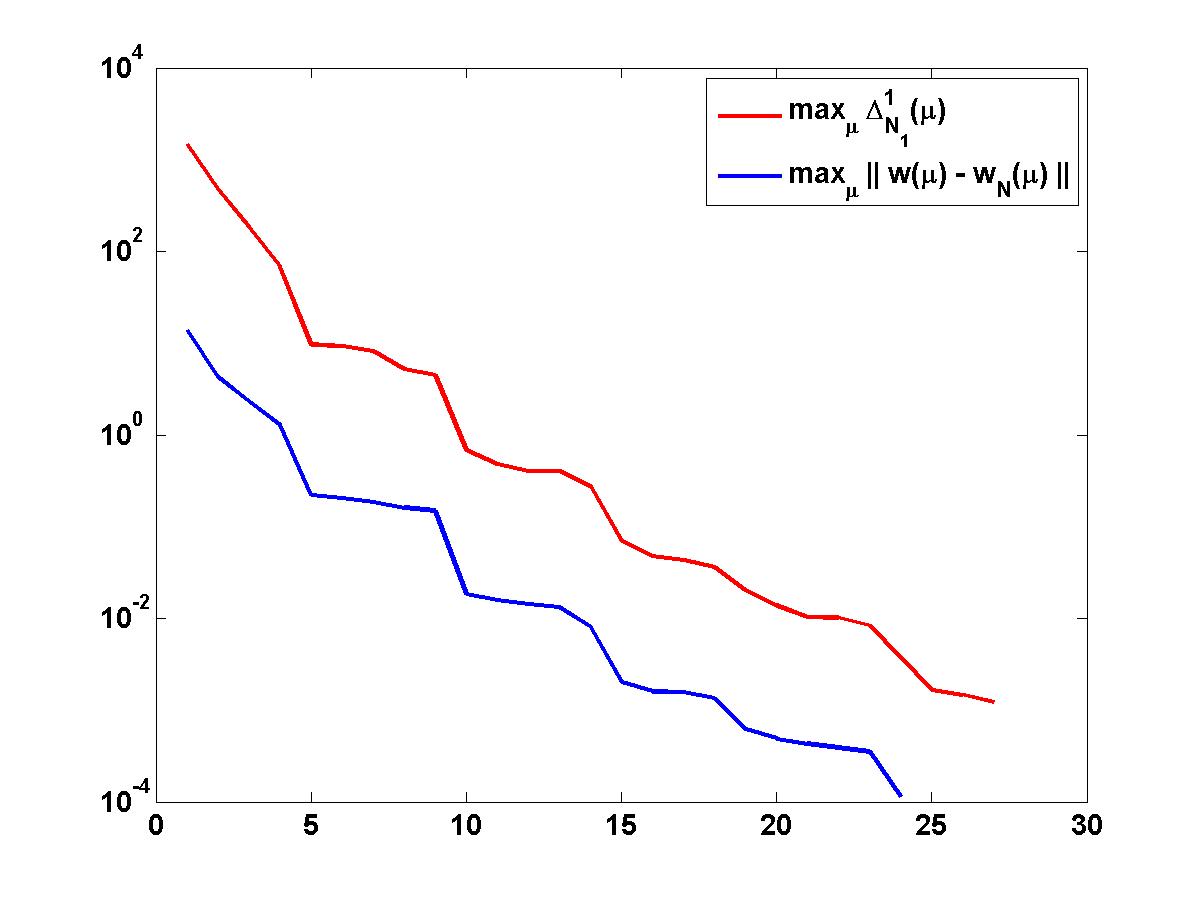}
				\label{fig:MaxErrSeq}}
	\hfill
	\subfigure[Error estimator and true error for $\rho \in (-1,1)$]{
		\includegraphics[width=0.45\textwidth]{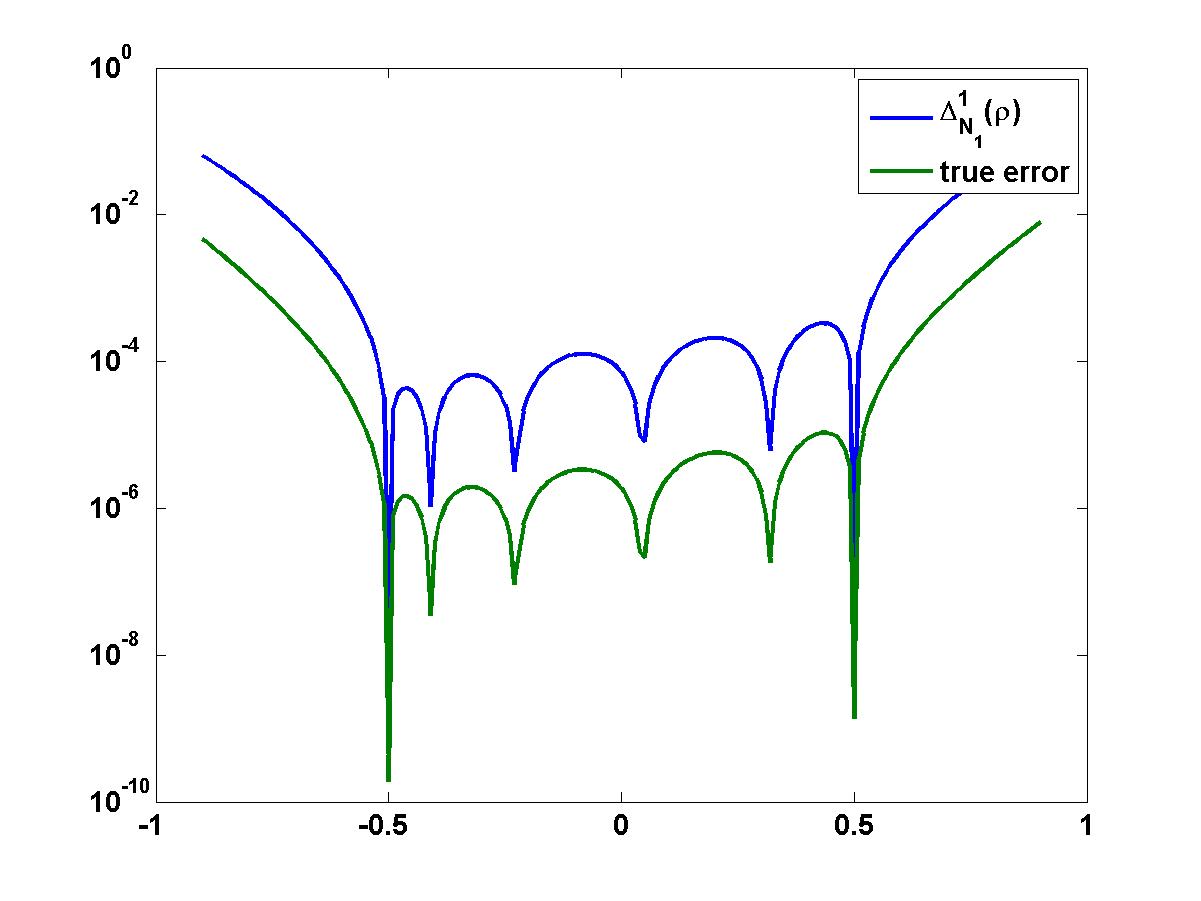}
		\label{fig:Rho}}
	\caption{Performance of error estimators\label{Fig:ErrEst}}
\end{figure}

We compare the evolution greedy with the proposed error estimate $R_{N_1,1}$ with the so called strong (evolution) greedy using the exact error (instead of the estimate) which is determined using a detailed solution. Doing so, we can investigate the performance of the error estimator. As we see in Figure \ref{Fig:ErrEst} \subref{fig:MaxErrSeq}, the error estimate behaves similar to the true error.  An error tolerance of $10^{-3}$ is reached by $28$ basis functions using the estimate, whereas the strong greedy shows that $24$ basis functions suffice to reach a tolerance of $10^{-4}$. As expected, the error bound overestimates the error. One reason is that (for simplicity) we use a pessimistic lower bound of the inf-sup condition ($\beta_{\text{LB}} = 0.005$).

\subsubsection{Error propagation}
Now, we investigate how different training sets influence the results. 
First, we choose $u_0 = h^1$ (i.e., an exact initial condition). As we see in Figure \ref{Fig:ErrEst} \subref{fig:Rho} the error for $\rho \in [-0.5,0.5]$ (the part that is covered by the training set) is below $10^{-4}$ and the error estimation $R_{28,1}$ is below $10^{-3}$. For $\rho \in (-1,1) \setminus [-0.5,0.5]$ the approximation is --as expected-- worse. Note that we ploted on a grid for $\rho$ with step size $0.01$. Furthermore, we can see that $(h^1, -0.5)$, $(h^1, -0.4091)$, $(h^1, -0.2273)$, $(h^1, 0.0455)$,  $(h^1, 0.3182)$ and $(h^1, 0.5)$ are in the sample set for $\W_{28}$.

\begin{figure}[!htb]
	\subfigure[Detailed $u^\cN$ at $t = 0$.]{\label{fig:truthInit}
	  				\includegraphics[width = 0.48\textwidth]{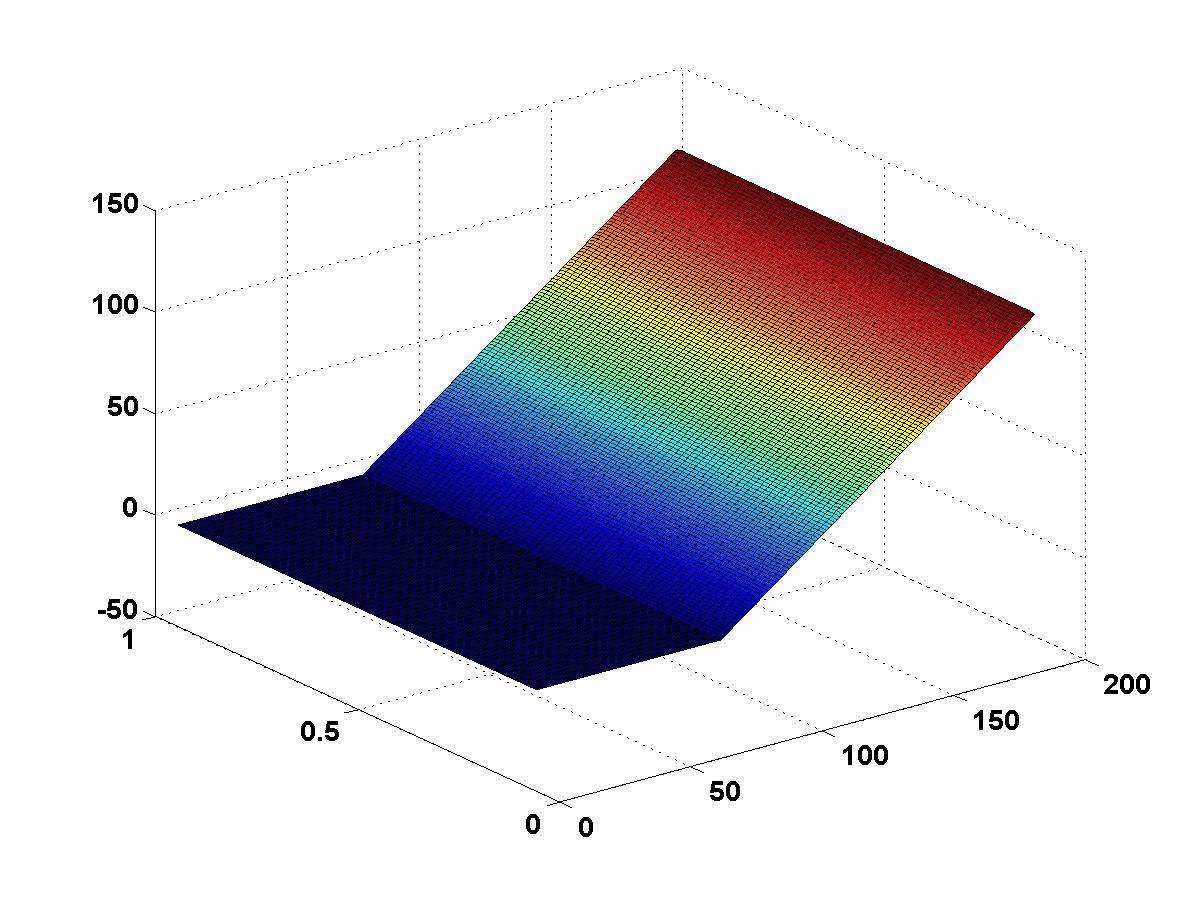}}	
	\subfigure[Reduced $u_N$ at $t = 0$.]{\label{fig:RBInit}
	  				\includegraphics[width = 0.48\textwidth]{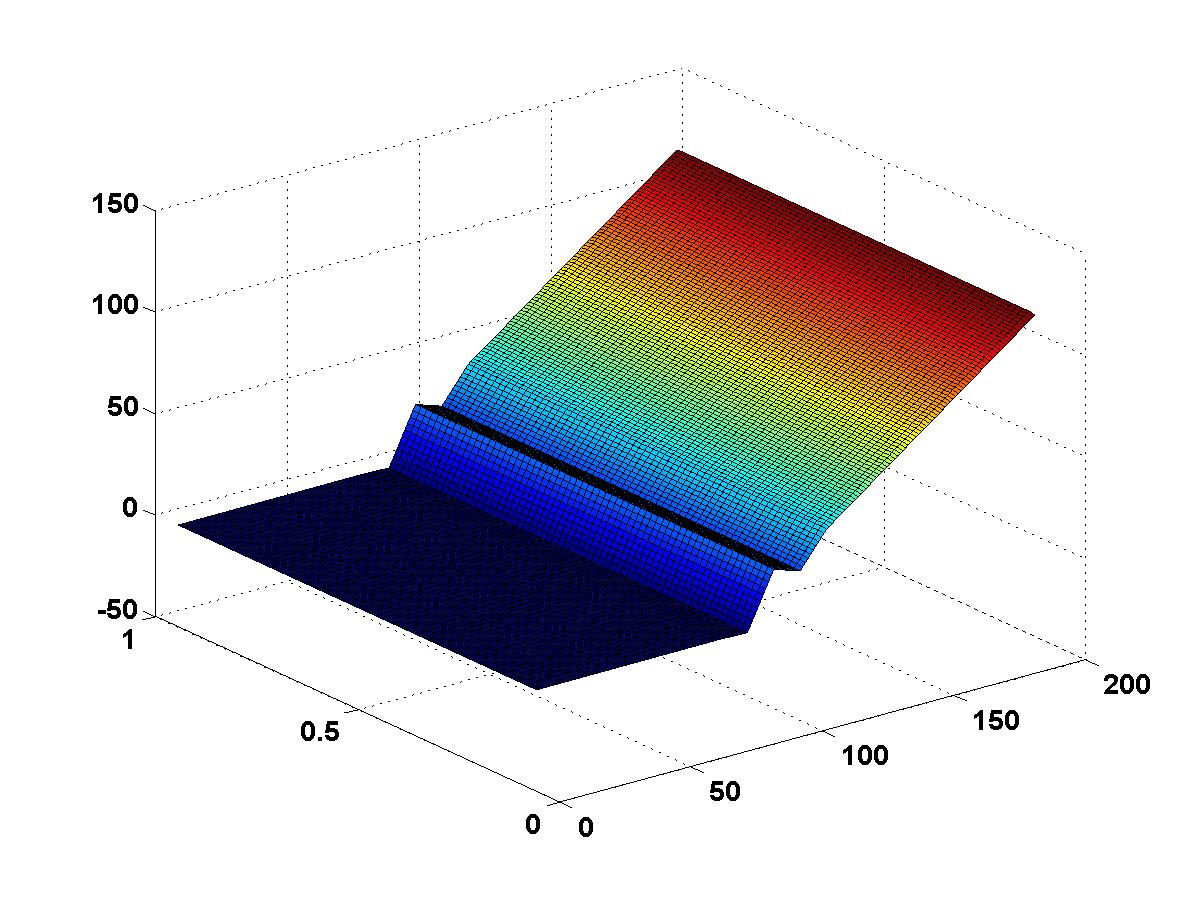}}
	\caption{\label{Fig:u0}Detailed and reduced initial condition ($u |_{[10^{-8},190] \times [0.05,0.95]}$).}
\end{figure}

Next, we choose a call payoff function $\mu_0(y,\nu) = \max(\exp(y)-K,0)$ with strike price $K = 70$ and $\rho = 0.3$. The resulting approximation error using  only  five POD basis functions can clearly be seen in Figure \ref{Fig:u0}. Obviously, the space-time error is large on the domain $\Omega$, in particular due to the errors near the boundaries. However, usually one has to enlarge the domain in advance due to the chosen boundary conditions, hence one is only interested in a smaller part of $\Omega$. 

\begin{figure}[here]
			\includegraphics[width=0.45\textwidth]{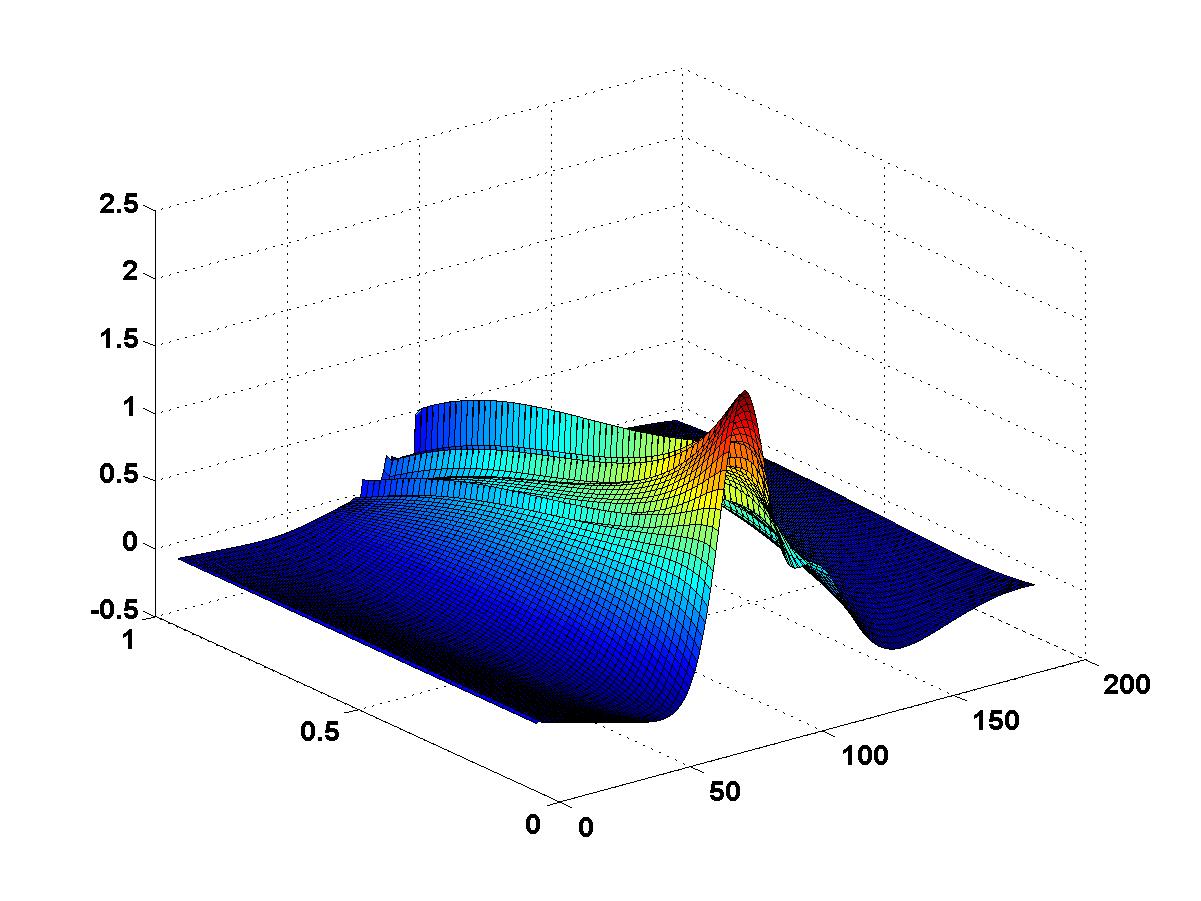}
			\caption{Absolute error $u^\cN (T) - u_N (T)$ at maturity ($u |_{[10^{-8},190] \times [0.05,0.95]}$).}
			\label{fig:absError}
\end{figure}

For option pricing the \emph{absolute error}, i.e., the difference of the prices, is also of interest, because the resulting RB price is to be actually paid by the customer. As we can see in Figure \ref{fig:absError}, the absolute error at the final time $T$ is actually about $2$, which is clearly too large. The reason is the POD truncation of the initial value based upon a \emph{relative error} basis. 

\begin{figure}[!htb]
	\subfigure[Maximum errors over iterations using the error estimate $R^1_{N_1}$.]{\label{fig:MaxErrVgl}
	  				\includegraphics[width = 0.48\textwidth]{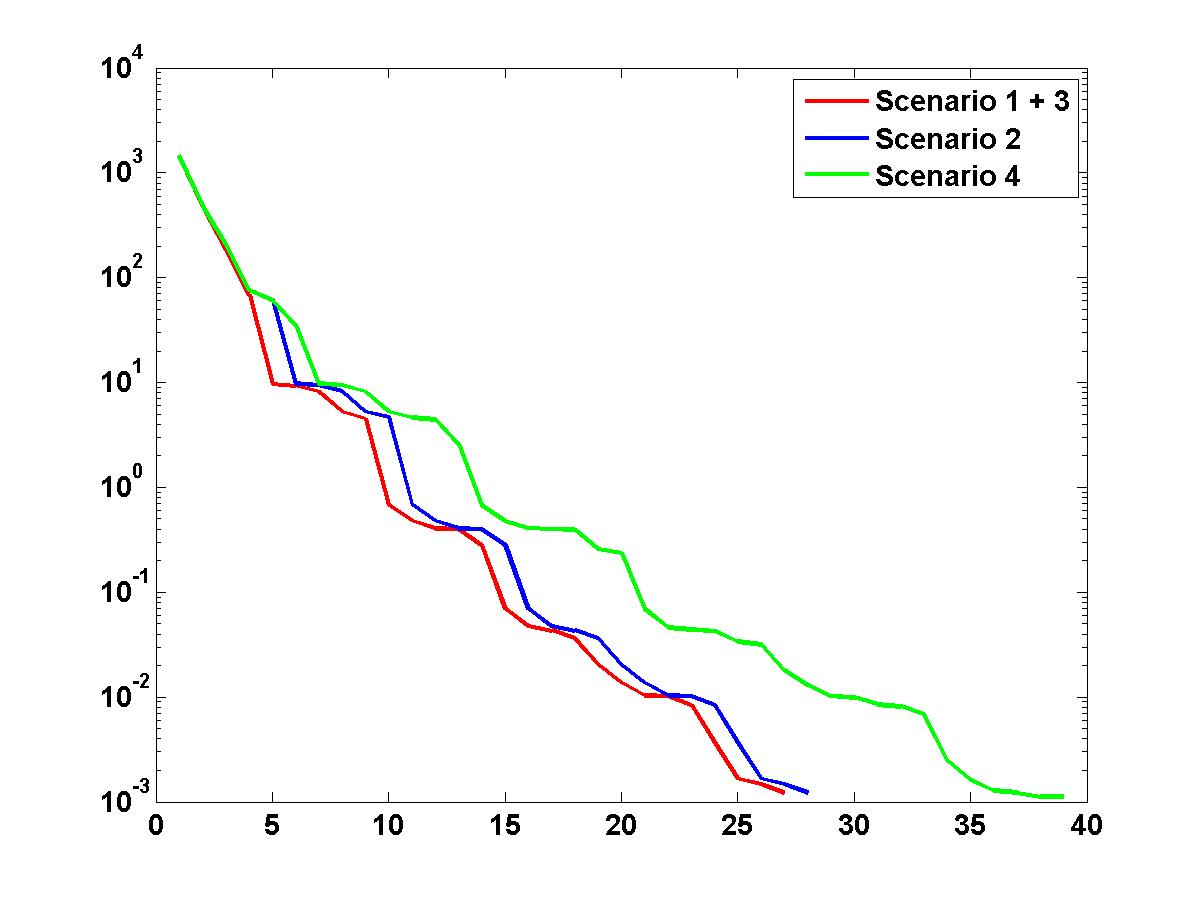}}	
	\subfigure[True errors for $\rho \in \lbrack-0.5,0.5 \rbrack$.]{\label{fig:TrueErrVgl}
	  				\includegraphics[width = 0.48\textwidth]{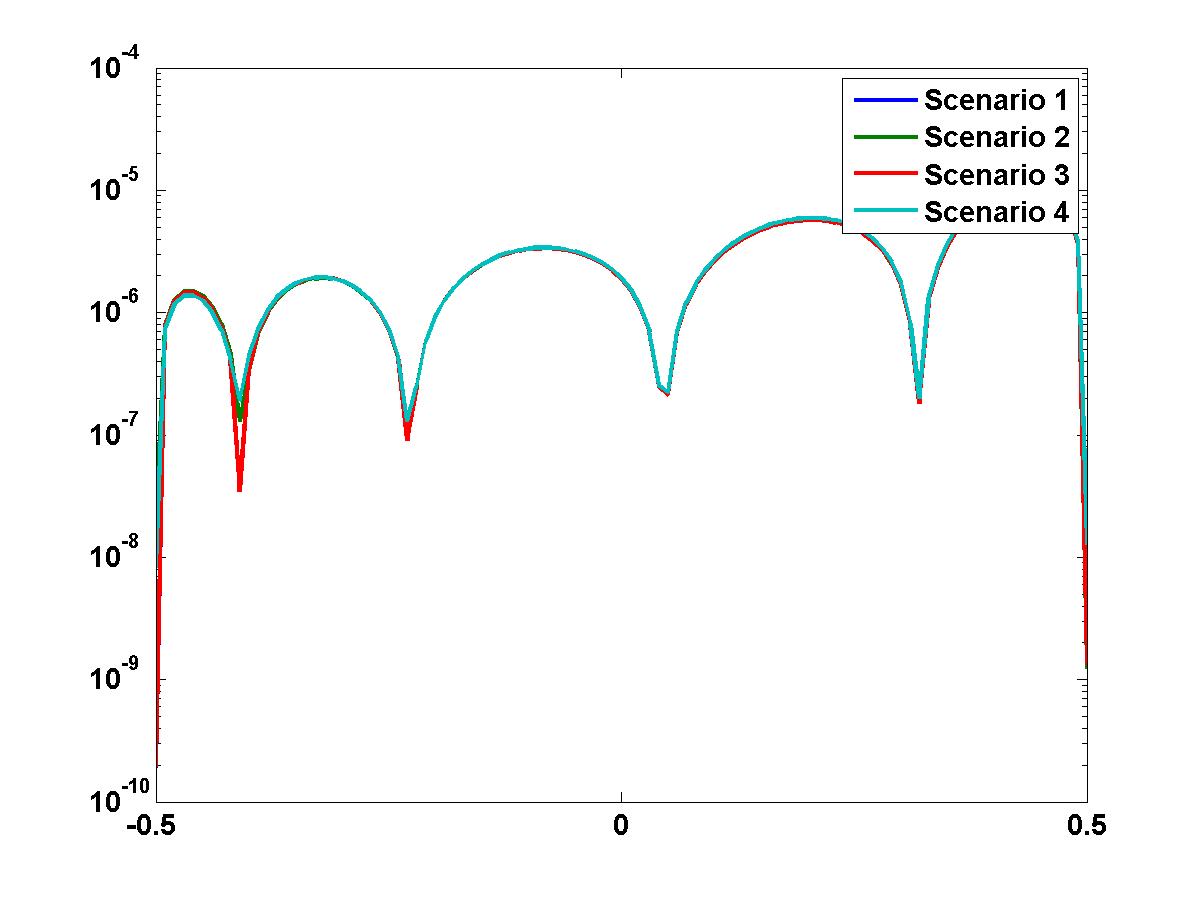}}
	\caption{Comparison of the different scenarios presented in Table \ref{tab:Mtrain}.}
\end{figure}

As mentioned before, we can choose the training set $M_{\text{train}}$ in the evolution greedy independend of $S_{N_0}^0$. We compare different combinations of reduced basis spaces in Table \ref{tab:Mtrain}. For $\rho \in [-0.5,0.5]$ with $u_0 = h^1$ (so again an exact initial condition) the true errors of the different approximations are shown in Figure \ref{fig:TrueErrVgl}. We can see, that for the exact initial value all RB approximations are evenly good. In contrast, by looking at the last column of Table \ref{tab:Mtrain}, we note that if we choose a poor initial value approximation resp.\ a poor training set for $\cD_0$, the RB approximation is not acceptable. In particular for PDEs with zero right-hand side and a linear operator (as we have for option pricing), we can show that for every new parameter $\bar{\mu}$, where the function part $\bar{\mu}_0$ lies in the span of the training set $M_{\text{train}}^1$, we get $\Delta_{N_1}^1(\bar{\mu}) < \text{tol}_1$. That explains why we do not get a larger $N_1$ resp. a better RB approximation by extending $M_{\text{train}}^1$ in Scenario 2 in comparison to Scenario 1. Of course, in this situation, this is also explained by  Gronwall's lemma. 

 \begin{table}[!htb]
 \caption{Comparison of different sets for $H_{N_0}$ and the training set of the evolution greedy ($M_{\text{train}} = M_{\text{train}}^0 \times M_{\text{train}}^1$, $M_{\text{train}}^1$ as before, $\text{tol}_1 = 10^{-3}$, $\mu(y,\nu) = (\max ( e^y -K,0), 0.3)$ and with $u|_{[10^{-8},190] \times [0.05,0.95]}$). }
 \begin{tabular}{c|c|c|c|c|c}
   Scenario & $H_{N_0}$ & $M_{\text{train}}^0 $ & $\# M_{\text{train}}$ & $N_1$ & $\| u^\cN(\mu) - u_N(\mu)\|_{\bar{\X}}$\\
	 \hline 1 & $\Span \{ h^1, \ldots, h^5 \}$ & $\{ \mu_0^1, \ldots, \mu_0^5 \}$ & $60$ & $28$ & $9.1994$ \\
	 2 & $\Span \{ h^1, \ldots, h^5 \}$ & $\{ \mu_0^1, \ldots, \mu_0^7 \}$ & $84$ & $29$ & $9.1994$ \\
	 3 & $\Span \{ h^1, \ldots, h^7 \}$ & $\{ \mu_0^1, \ldots, \mu_0^5 \}$ & $60$ & $28$ & $7.9801$ \\
   4 & $\Span \{ h^1, \ldots, h^7 \}$ & $\{ \mu_0^1, \ldots, \mu_0^7 \}$ & $84$ & $40$ & $1.7162e-04$ \\
 \end{tabular}
 \label{tab:Mtrain}
 \end{table}


\mysection{Summary and outlook}
\label{Sec:6}
We have introduced a Reduced Basis space-time variational approach for parametric parabolic partial differential equations having coefficient parameters and a variable initial condition. Feasibility and efficiency have been demonstrated. Obviously, a whole variety of further questions arises, just to mention some of them that we aim to consider in the near future: \\[-4mm]
\begin{itemize}
	\item Extension of the Bernstein representation of the initial condition to an adaptive wavelet approximation.
	\item In \cite{CS}, an alternative space-time variational formulation has been considered which transfers the essential initial condition to a natural one.
	\item In case of a basket of options, the dimension of the problem of course grows, which calls for a specific treatment e.g.\ by the Hierarchical Tucker format as e.g.\ in \cite{DK,KRU}. More general, the choice of trial and test spaces within the (offline) phase can be further investigated and possibly optimized.
	\item For calibration purposes, all model parameters have to be taken into account. This type of high dimensionality needs particular treatment as e.g. in \cite{DMR}, \cite{HSZ} or \cite{HZ}. Moreover, numerical stabilization techniques have to be investigated.
	\item Extensions to other financial models, American options, etc.
\end{itemize}


\end{document}

%% file: KU-0header.tex
\usepackage[latin1]{inputenc}
\usepackage[english]{babel}
\usepackage{exscale}
\usepackage{dsfont}
\usepackage[plain]{algorithm}
\usepackage{algorithmic}
\usepackage{amsmath, amssymb}
\usepackage{amsthm}
\usepackage{mathrsfs}
\usepackage{latexsym}
\usepackage{wasysym}
\usepackage{framed}
\usepackage{color}
\usepackage{ifpdf}
\usepackage{rotating}
\usepackage{subfigure}
\usepackage{multirow}
\usepackage{multicol}

\usepackage{paralist}

\usepackage[ulem=normalem,final]{changes}

\theoremstyle{plain}
\newtheorem{thm}{Theorem}[section]

\newtheorem{prop}[thm]{Proposition}

\newtheorem{rem}[thm]{Remark}

\theoremstyle{definition}

\theoremstyle{definition}

\floatstyle{ruled}
\newfloat{algo}{h!}{loa}
\floatname{algo}{Algorithm}

\newcommand{\cO}{\mathcal{O}}

\newcommand{\R}{\mathbb{R}} 
\newcommand{\N}{\mathbb{N}} 
\newcommand{\W}{\mathbb{W}} 
\newcommand{\X}{\mathbb{X}} 
\newcommand{\Y}{\mathbb{Y}} 
\newcommand{\Z}{\mathbb{Z}} 
\newcommand{\Q}{\mathbb{Q}} 

\numberwithin{equation}{section}

\renewcommand{\sectionmark}[1]{}

\newcommand{\BIGOP}[1]{\mathop{\mathchoice%
{\raise-0.22em\hbox{\huge $#1$}}%
{\raise-0.05em\hbox{\Large $#1$}}{\hbox{\large $#1$}}{#1}}}

\newcommand{\mysection}[1]{\section{#1}
\setcounter{subsubsection}{0}
\setcounter{equation}{0}
\setcounter{thm}{0}
\setcounter{table}{0}
\setcounter{figure}{0}
}

\newcommand{\cD}{\mathcal{D}}

\newcommand{\cI}{\mathcal{I}}
\newcommand{\cJ}{\mathcal{J}}
\newcommand{\cK}{\mathcal{K}}
\newcommand{\cL}{\mathcal{L}}
\newcommand{\cM}{\mathcal{M}}
\newcommand{\cN}{\mathcal{N}}

\renewenvironment{itemize}{\vspace*{-13pt}
  \begin{list}{$\bullet$}
  { \setlength{\itemsep}{0pt}
     \setlength{\parsep}{1.5pt}
     \setlength{\topsep}{0pt}
     \setlength{\partopsep}{0pt}
     \setlength{\leftmargin}{1.2em}
     \setlength{\labelwidth}{1.2em}
     \setlength{\labelsep}{0.5em} 
     }\item[]}{\vspace*{0mm}\end{list}}

\usepackage{stmaryrd}
\usepackage{url}
\usepackage{mathtools}
\usepackage[normalem]{ulem}

\usepackage{algorithmic}

\usepackage[colorinlistoftodos]{todonotes}

\newcommand\Init{{\rm{init}}}
\newcommand\Space{{\rm{space}}}
\newcommand\Time{{\rm{time}}}
\newcommand\Span{{\rm{span}}}
\newcommand\Div{{\rm{div}}}

\newcommand{\Dt}{{\Delta t}}

\newcommand{\bbf}{\mathbf{f}}
\newcommand{\bg}{\mathbf{g}}
\newcommand{\bh}{\mathbf{h}}
\newcommand{\bq}{\mathbf{q}}
\newcommand{\bs}{\mathbf{s}}
\newcommand{\bu}{\mathbf{u}}

\newcommand{\bw}{\mathbf{w}}
\newcommand{\bz}{\mathbf{z}}

\newcommand{\bmu}{\boldsymbol{\mu}}

\newcommand{\bA}{\mathbf{A}}
\newcommand{\bB}{\mathbf{B}}
\newcommand{\bC}{\mathbf{C}}

\newcommand{\bM}{\mathbf{M}}
\newcommand{\bN}{\mathbf{N}}

\usepackage{color}

\usepackage[english]{babel}
\usepackage{stmaryrd}

\newcommand\skp[2]{(\kern-2pt(#1,#2)\kern-2pt)_X}

\definecolor{lred}{RGB}{231,91,18}
\definecolor{dred}{RGB}{210,0,0}
\definecolor{dgreen}{RGB}{0,135,0}
\definecolor{lblue}{RGB}{0,120,225}
\definecolor{grey}{RGB}{175,175,175}
\definecolor{orange}{RGB}{250,140,0}
\definecolor{cyan}{RGB}{0,215,215}
\definecolor{white}{RGB}{255,255,255}